\documentclass[preprint,9pt]{elsarticle}
\usepackage[fleqn]{amsmath}

\usepackage{amssymb}
\usepackage{amssymb} 

\usepackage{amsmath}
 \usepackage{amsthm}
 \usepackage{hyperref}
 \usepackage{natbib}

\usepackage{graphicx}

\usepackage[all]{xy}

\newcommand{\N}{\mathbb{N}}
\newcommand{\R}{\mathbb{R}}

\newcommand{\C}{\mathbb{C}}
\newcommand{\p}{\mathbf{P}}

\def\norm#1{\left \Vert#1\right \Vert}

\def\abs#1{{\vert #1\vert}}

\def\text#1{\mbox{#1}}

\newtheorem{lemma}{ \kern \parindent Lemma}
\newtheorem{definition}{ \kern \parindent Def{}inition}
\newtheorem{proposition}{ \kern \parindent Proposition}
\newtheorem{corollary}{ \kern \parindent Corollary}
\newtheorem{remark}{ \kern \parindent Remark}
\newtheorem{example}{ \kern \parindent Example}

\begin{document}

\begin{frontmatter}

\title{Algorithms for computing basins of attraction associated with a rational self-map of the Hopf fibration based on Lyapunov exponents\tnoteref{mytitlenote}}
\tnotetext[mytitlenote]{Partially supported by University of La Rioja and project: PID2020-118753GB-I00 of the Spanish  Ministerio de Ciencia e Innovaci\'on. }

\author{V. \'Alvarez-Aparicio\fnref{label1}}
\ead{valvarezaparicio@gmail.com}
\author{J. M. Garc\'{\i}a-Calcines\fnref{label2}}
\ead{jmgarcal@ull.edu.es}
\author{L. J. Hern\'andez-Paricio\corref{cor1}\fnref{label1}}
\cortext[cor1]{Corresponding author}
\ead{luis-javier.hernandez@unirioja.es}
\author{M. T. Rivas-Rodr\'{\i}guez\fnref{label1}}
\ead{maria-teresa.rivas@unirioja.es}


\affiliation[label1]{organization={Universidad de La Rioja},
addressline={Dpto. de Matem\'aticas y Computaci\'on, Edificio CCT, C/Madre de Dios, 53},
city={Logro\~no},
postcode={26006},
country={Spain}
}

\affiliation[label2]{organization={Universidad de La Laguna},
 addressline={Dpto. de Matem\'aticas, Estad\'{\i}stica e I.O.  e Instituto de Matem\'aticas y Aplicaciones (IMAULL), Avda. Astrof\'{\i}sico Fco. S\'anchez, s/n., Facultad de
Ciencias. Sección de Matem\'aticas},
city={La Laguna},
postcode={38200},
country={Spain}
}

\medskip




\begin{abstract}

{\it Dedicated to Jos\'e Manuel Rodr\'{\i}guez Sanjurjo on the occasion of his 70th birthday.}

\medskip

The objective of this work is twofold. One, we develop a theory of iteration for complex  rational functions so that the problems of overflow and indeterminacy caused by null, or almost null, denominators can be avoided when developing implementations. Second, we present easily implementable methods that allow the calculation of attracting cycles  as well as the graphical representation of their basins of attraction.

In order to deal with our first goal we work with homogeneous complex coordinates and we take the complex projective line as a model, which is expressed as a quotient of the 3-sphere through the Hopf fibration. An irreducible representation of a rational function can now be presented as a Hopf fibration endomorphism.
As well as our second goal is concerned, we use Lyapunov functions and exponents to calculate the cycles associated with  endomorphisms and to graphically represent the corresponding basins of attraction. We point out that our algorithms are based on the calculation of a finite set of non-negative real constants and their calculation does not depend on the previous determination of the fixed points or attracting cycles.

\bigskip

\end{abstract}




\end{frontmatter}



\section{Introduction}\label{In}

In  this work, we have two main objectives:
\begin{itemize}
\item Avoid the overflow (and underflow) errors that the computational iteration of a rational function
presents when either the denominator is close to zero or when both the numerator and denominator are very small.
  \item Give a graphical representation of the basins of attracting cycles by using Lyapunov functions and exponents.
\end{itemize}

Let $\C[x]$ denote the ring of polynomials with complex coefficients and $\C(x)$ the field of complex rational functions
$f(x)=\frac{A(x)}{B(x)}$, where $A, B \in \C[x]$ and $B$ is not the zero polynomial. When one tries to find the values of the expression $\frac{A(x_0)}{B(x_0)}$ in a computational system, there may be two different kind of issues:

\bigskip

a) Mathematical difficulties:

\begin{itemize}
  \item[(i)] If $A(x_0)\not = 0$ and $B(x_0) = 0$, then $\frac{A(x_0)}{B(x_0)}$ is not defined ($x_0$ is a pole of the function $f$).
  \item [(ii)] If $A(x_0) = 0$ and $B(x_0) = 0$, then $\frac{A(x_0)}{B(x_0)}$ is not defined ($x_0$ is an indetermination point of $f$).
\end{itemize}

b) Computacional  difficulties:

\begin{itemize}

\item[(iii)] If $A(x_0)\not = 0$ and $\abs{B(x_0)}$ is very small (that is, $\abs{B(x_0)}\leq 10^{-N}$ for a large natural number $N$), then an overflow error can be produced when the computational system tries to compute $\frac{A(x_0)}{B(x_0)}$.
\item [(iv)] If $\abs{A(x_0)}$ or $\abs{B(x_0)}$ are very large, then an overflow error may occur.

\item [(v)] If $\abs{A(x_0)}$ and $\abs{B(x_0)}$ are very small (that is, $\abs{A(x_0)}+\abs{B(x_0)}\leq 10^{-N}$ for a large natural $N$), then an indetermination error can be produced when trying to compute $\frac{A(x_0)}{B(x_0)}$.
\end{itemize}

The mathematical difficulties can be solved by using suitable models of spaces and functions:

For the purpose of solving problem (i), instead of working with $\C$ we can use the Alexandroff compactification of $\C$, denoted as $\hat \C = \C \cup \{\infty\}$, and other equivalent models such as the unit 2-sphere $S^2$ or the complex projective line $\p^{1}(\C)$ (see subsection \ref{structCinf}). The Euclidean metric of $S^2$ gives rise to a metric space structure to each of these models. Moreover, the structure of 1-complex manifold in $\p^{1}(\C)$ permits us to consider the irreducible rational maps as analytic maps.  We observe that if we take homogeneous coordinates $x=\frac{z}{t}$ and $d=\max\{\mbox{degree} (A(x)), \mbox{degree} (B(x))\}$, then

$$\frac{A(x)}{B(x)}=\frac{t^d A(z/t)}{t^dB(z/t)}.$$

We also have that $F(z,t)=t^d A(z/t)$ and $G(z,t)=t^dB(z/t)$ are homogeneous  bivariate polynomials of degree $d$ verifying

$$\frac{A(x)}{B(x)}=\frac{F(z,t)}{G(z,t)}.$$

In this case we obtain a pair $R=(F,G)$ where $F(z,t),G(z,t)$ are homogeneuos bivariate polynomial with the same degree. It is important to remark that  we have a map $R\colon \C^2 \to \C^2$.  Taking into account that  $R$ is a homogeneous pair (see Definition \ref{homogeneouspair}), when $R$ is $r$-irreducible (see Definition \ref{irreducible}), we also have an induced map $R^P\colon \p^1(\C) \to \p^1(\C)$ (see definition given by equation \eqref{pe}).

\bigskip

Therefore problem (i) can be bypassed by taking into account that:

\begin{itemize}

\item[(i)]
 The use of homogeneous coordinates allows us to work at the point at infinity; the representation of a rational function by a pair of homogeneous  bivariate polynomials with the same degree allows us to compute the numerical value of the function at any pole point and at the infinity point.
\end{itemize}

With regard to problem (ii), the mathematical model denominated {complex projective line} $\p^{1}(\C)$ can  be slightly altered with the addition of a new element: { the super-zero point} $\bar 0$.  In this way, we get what we call the {augmented projective line} $\p^{1+}(\C)=\p^{1}(\C)\cup \{\bar 0\}$ as well as the corresponding { augmented 2-sphere} $S^{2+}=S^2\cup \{(0,0,0)\}$.  The new models are quite compatible when a rational map is represented by a pair of homogeneous bivariate polynomials with the same degree.
 \begin{itemize}
\item[(ii)]
This permits us to solve the indetermination problem (ii) which appears when the values of the numerator and denominator of a rational function are simultaneously equal to zero. In these cases, the super-zero point is taken as the image of an indeterminate point.
\end{itemize}

Using the augmented models $\p^{1+}(\C)$ and $S^{2+},$ a non irreducible representation of a rational map induces self-maps of the augmented models. The iteration of these augmented maps are closely related to the iterations of the corresponding irreducible representation of the rational map. For a study of the use and applications of augmented models we refer the reader to \cite{GGHR15}, \cite{HernndezParicio2015BivariateNM} and \cite{EGHR2021}.

Nevertheless, in this work we are avoiding problem (ii) by just working with irreducible representations of rational functions.
Observe that if $f(x)=\frac{A(x)}{B(x)}$ is a rational map and $x_0$ is a common zero of $A(x)$ and  $B(x)$, then $A(x)=(x-x_0)A_1 (x)$ and $B(x)=(x-x_0)B_1(x)$.  Iterating  this procedure we can find an irreducible representation $g(x)=\frac{\tilde A(x)}{\tilde B(x)}$ without indetermination points. Therefore, $\frac{\tilde A(x)}{\tilde B(x)}$ is a new representation of the rational function avoiding indetermination problems. Also notice that the map $g$ only differs from the map $f$ on the finite set of indetermination points of $g$.

\begin{itemize}

\item[(iii)] The  computational  problems  (iii) disappear if we use homogeneous coordinates $[z:w]$. Indeed, observe that
$[A(x_0):B(x_0)]=[F(x_0,1):G(x_0,1)]$ and if $B(x_0)$ is very small, then the computational system changes from $[F(x_0,1):G(x_0,1)]$ to $[F(x_0,1):0)]$, which represents the point at infinity in $\p^{1}(\C)$.

\item[ (iv)] Problem (iv) can be eluded by using normalized coordinates (see Definition \ref{normalizedhc}) and considering the normalized representation of a homogenous pair (see subsection \ref{normphi}). Indeed, for every pair $(z,t)$ satisfying $\abs{z}+\abs{t}=1$ we have that
$[A(z/t):B(z/t)]=[F(z,t):G(z,t)]$ verifies the inequality $\abs{F(z,t)}+\abs{G(z,t)}\leq \norm{(F,G)}$. Furthermore, if we take a normalized representation $( {\bar F}, {\bar G})$ of $\frac{A(x)}{B(x)}=\frac{A^*(x)}{\bar B^*(x)}$, then
$[A(z/t):B(z/t)]=[ A^*(z/t): B^*(z/t)]=[\bar F(z,t): \bar G(z,t)]$ and
$\abs{\bar F(z,t)}+\abs{\bar G(z,t)}\leq \norm{(\bar F,\bar G)}=1.$

Therefore, we deal with the computation of $(\bar F(z_0,t_0), \bar G(z_0,t_0))$ ($x_0=z_0/t_0$) where $\abs{\bar F(z_0,t_0)}+\abs{\bar G(z_0,t_0)} \leq 1$ rather than trying to compute $A(x_0)/B(x_0)$, where $A(x_0)$ and $B(x_0)$ can be large complex numbers.

\item[(v)] If we take a normalized representation  $(\bar F, \bar G)$ for an irreducible representation of a rational map, then
we obtain a continuous map $ S^3 \to [0,1]$ defined as
$(z,t) \mapsto \abs{\bar F(z,t)}+\abs{\bar G(z,t)}$.
Since $S^3$ is a compact space and $\abs{\bar F(z,t)}+\abs{\bar G(z,t)}\not =0$, for every $(z,t) \in S^3$, there exists $N\in \N$ such that
$10^{-N}  \leq \abs{\bar F(z,t)}+\abs{\bar G(z,t)}$, for all $(z,t) \in S^3$.

Consequently, in a suitable computational environment, which is able to work up to precision of order $10^{-N }$, we can compute the normalization of ${(\bar F(z,t)},{\bar G(z,t))}$ to obtain again a point in $S^3$.

\end{itemize}

Therefore, given  an irreducible quotient of two polynomials  we can consider the maps $ R^S, R^P$ (see formulas \eqref{ese}, \eqref{pe}) that give a self-map of the Hopf fibration:

$$\xymatrix{
\ar[d]_q S^{3} \ar[r]^{R^S}  & S^{3} \ar[d]^q \\
\p^{1}(\C) \ar[r]^{R^P} &\p^{1}(\C). }
$$

Our methods and algorithms are based on the study of the iteration of this self-map $( R^S, R^P)$ of the Hopf fibration.

\medskip

The second main objective of this paper is to compute the basins of attraction  of a rational function $f$. In the literature there are many algorithms whose first part involves the calculus of the finite set of attracting cycles (see, for instance, \cite{Var02}, \cite{gutierrez_hernandez_magrenan_rivas_2016}). A second part of these algorithms takes a point and one study whether its orbit converges or not to some of these attracting cycles.

In this paper we consider a different method, which depends on the choice of a suitable continuous map $\phi$ from the Riemann sphere to $\mathbb{R}_+=[0,+\infty )$, the space of non-negative real numbers. 
Then we  consider what we call the Lyapunov function
$L_f(\phi)$ induced by $f$ and $\phi$, see Definition \ref{Lyafunction}. In general, $L_f(\phi)$ is a non-continuous function which is constant almost everywhere on each basin of an attracting $n$-cycle. The first part of our algorithm is devoted to computing the finite set of Lyapunov constants: 
$$0 \leq c_1< c_2 < \cdots < c_{l_-1}< c_l <  1\leq c_{l+1}< c_{l+2}<\cdots<c_{l+k-1}< c_{l+k}.$$ 

In the second part it computes, up to some tolerance, the value $L_f(\phi)(x)$, of any point $x$ where $L_f(\phi)$ is defined. The algorithm also finds, when it is workable, the position $i \in \{1, \cdots, l\}$ such that $L_f(\phi)$ agrees (again, up to some tolerance) to $c_{i}$. In this way the Riemann sphere can be divided into $l+1$ disjoint regions, where each of which corresponds to one of these constants less than 1 and its complement. It is noteworthy to mention that, when the number of attracting cycles is equal to the number of constants which value is  less than 1, these regions are exactly the basins of attraction  of the attracting cycles. Moreover, under some suitable additional conditions, the complement of the union of the basins of attraction  is precisely the Julia set, from which one can find an ergodic measure satisfying that $L_f(\phi)$ is also constant almost everywhere.

In principle, one may consider different viable choices for a continuous map $\phi \colon \p^1(\C) \to \R_+$. However, for our study of the iteration of a rational map $f\colon \p^1(\C) \to \p^1(\C)$ we have chosen the \emph{spherical derivative} $\phi=f^{\#} \colon \p^1(\C) \to \R_+$ (for a study of a different choice of $\phi$, we refer the reader to \cite{EGHR2021}).
We note that if we take the logarithms of the constants induced by the spherical derivative:
$$-\infty \leq \log c_1< \cdots < \log c_{l_-1}< \log  c_l <  0 \leq \log c_{l+1}<\cdots<\log c_{l+k-1}< \log c_{l+k},$$ 
 we obtain the Lyapunov exponents of the rational map $f$, where negative, zero or positive values correspond  to attracting, indifferent or repelling cycles, respectively. For more results on Lyapunov exponents for rational maps, see \cite{GPR10}, \cite{DeMarco03}.
 
 \medskip

In summary, we have developed a procedure that avoids the indetermination problems, which may  appear computing the spherical derivative, by using a suitable pair of bivariate functions, see Proposition \ref{computinspherical}.  We also have developed a set of algorithms (see section \ref{algorithms}) that permits us to compute the basins of attraction  of cycles by using the Lyapunov function associated to the spherical derivative. As an example of application of the implemented algorithms developed in this paper, we have included an study of the Newton's method associated to cubic polynomials with simple roots, see section \ref{cubic_polynomials}, and we have taken an example of rational Hopf endomorphim giving the attraction basin of a 2-cycle (see subsection \ref{example2basin}). 

The procedure and algorithms presented in this paper have been implemented in Julia (see \cite{Alvarez-Aparicio_Lyapunov_Cycle_Detector_2022}) following ideas of a previous implementations in Mathematica (see \cite{ExtremianaAldana2018SelfoverlaysAS}) and in early (deprecated) Julia versions (see \cite{hernandez19}).
Julia  programming language has automatic translation of formulas into efficient executable codes. It produces fast, low-level machine code that has traditionally  been
generated only by static languages 
 (see \cite{BEKS17}); as a consequence,  the execution times are very fast when compared to other  implementations in other programming languages, that have been developed by the authors 
(see \cite{HRM2012}).
Moreover, Julia provides a multiprocessing environment  to allow programs to run on multiple processes in separate memory domains at once.

\section{Mathematical framework}\label{mf}

In order to establish a theoretical foundation that holds and justifies the correct construction of our algorithms we shall use the mathematical notions, results and techniques described in this section.

\subsection{Convergence of arithmetic and geometric means  of sequences}

In an effort to ease the reading of the paper we have included this subsection of auxiliary results about the convergence of arithmetic and geometrical means of real sequences that will be needed. As these basic results are easily available in the literature their proofs are omitted.

\begin{lemma}\label{arith0}
If $(r_k)_{k\in \N}  $ is a sequence in $\{-\infty\}\cup \R$ and $\lim_{k} r_k=r\in \{-\infty\}\cup \R$, then there exists $j_0 \in \N$ such that
$$\lim_{n\to +\infty} \left ( \frac{1}{n-j+1} \sum_{i=j}^{n} r_i \right ) =r$$ \noindent for every $j\geq j_0$. Moreover, if $r=-\infty$, then we can take $j_0=0$.
\end{lemma}

Observe that the exponential map $\text{exp}\colon \R \to (0,+\infty)$ and the logarithmic function $\log \colon (0,+\infty) \to \R$ induce a canonical homeomorphism
$$\log  \colon [0, +\infty) \to \{-\infty\}\cup \R$$ \noindent by taking $\log(0)=-\infty$. Using the properties of the logarithms we obtain as a consequence:

\begin{corollary}\label{prod1}
If $(s_k)_{k\in \N}  $ is a sequence in $[0,+\infty)$ and $\lim_{k} s_k=s \in [0,+\infty)$, then there exists $k_0 \in \N$ such that $$\lim_{n\to +\infty} \left (\prod_{i=k}^{n} s_i \right )^{ \frac{1}{n-k+1} }=s$$ \noindent for every $k\geq k_0$. Moreover, if $s=0$, then we can take $k_0=0$.
  \end{corollary}

\begin{lemma}\label{arith3}
Let $(r_n)_{n\in \N}$ be a sequence in $\{-\infty\}\cup \R$ and let $q\geq 1$ be an integer. For each integer $k$ with $0\leq k \leq q-1$ consider the subsequence $r_{n_j^k}= r_{q j+k}$, $j \in \N$. Then,
$\lim_{n} r_n=r$ if, and only if,  $\lim_{j} r_{q j+k}=r$,  for each $k$ with $0\leq k \leq q-1$.
\end{lemma}

\begin{lemma}\label{arith4} Let $q>1$ be an integer and $(u_j^0)_{j\in \N}, \cdots, (u_j^{q_-1})_{j\in \N}$ real sequences such that $\lim_j(u_j^k)=u_k\in \{-\infty\}\cup \R $. Define a new ``union'' sequence
$(\tilde u_i)_{i \in \N}$ by  $\tilde u_i=u_{j_i}^{k_i}$ where $j_i$ and $k_i$ are the unique nonnegative integers such that $0\leq k_i\leq q-1$ and $i=q j_i+k_i$.  Then there exists $j_0 \in \N$ satisfying that
$$\lim_{n\to +\infty} \left ( \frac{1}{n-j+1} \sum_{i=j}^{n} \tilde u_i \right ) =\frac{u_0+ \cdots + u_{q-1}}{q}$$
\noindent for every $j\geq j_0$. Furthermore, if there is $l \in \{0, \cdots, q-1\}$ with $u_l=-\infty$, then we can take $j_0=0$.
\end{lemma}

Again, using the logarithmic function and applying Lemma \ref{arith4} we obtain:

\begin{corollary}\label{arith5} Let $q>1$ be an integer and $(v_j^0)_{j\in \N}, \cdots, (v_j^{q-1})_{j\in \N}$ sequences in $[0,+\infty)$ such that $\lim_j(v_j^k)=v_k\in [0,+\infty)$. Define a new ``union'' sequence
$(\tilde v_i)_{i\in \N}$ by  $\tilde v_i=v_{j_i}^{k_i}$ where $j_i$ and $k_i$ are the unique integers satisfying $i=q j_i+k_i$ and $0\leq k_i\leq q-1$. Then there exists $j_0 \in \N$ such that,
$$\lim_{n\to +\infty} \left ( \prod_{i=j}^{n} \tilde v_i \right )^{ \frac{1}{n-j+1}} =(v_0 \cdots  v_{q-1})^{\frac{1}{q}}$$
\noindent for every $j\geq j_0.$
Moreover, if there is $l \in \{0, \cdots, q-1\}$ such that $v_l=0$, then we can take $j_0=0$.
\end{corollary}

 \subsection{Lyapunov, Time and Space Average Functions on discrete semi-flows}\label{LyaTimeSpaceAverage}

If $X$ is a (pseudo) metric space and $g\colon X\to X$ is a continuous map, then $(X,g)$ (or just $X$, for short) is said to be a \emph{discrete semi-flow}. Given an integer $n\geq 0$, $g^n$ will denote the $n$-th composition $g \circ\cdots \circ g$ and $g^0=id_X$. We will also consider the usual notions of \emph{fixed point}, \emph{periodic point} and \emph{$q$-cyclic point} of $g$. Note that if $x \in X$ is a $q$-cyclic point, then $x$ is a fixed point of $g^q$.

Now, if $\phi \colon X \to [0,+\infty)$ is a continuous map, then we can define the functions of Lyapunov type:
$$l_g^-(\phi)(x) = \liminf_{n\to +\infty} \left ( \phi(g^n(x))\right ), \quad l_g^+(\phi)(x)= \limsup_{n\to +\infty} \left ( \phi(g^n(x))\right ),$$

$$L_g^-(\phi)(x)= \liminf_{n\to +\infty} \Big( \prod_{k=0}^{n-1} \phi(g^k(x)) \Big )^{\frac{1}{n}}, \quad L_g^+(\phi)(x)= \limsup_{n \to +\infty}\Big(  \prod_{k=0}^{n-1} \phi(g^k(x)) \Big )^{\frac{1}{n}}.$$

\bigskip
And for any given integer $q\geq 1$, we can also take
$$l_g^{q-}(\phi)(x)= \liminf_{n\to +\infty} \Big( \prod_{k=n-q}^{n-1} \phi(g^k(x)) \Big )^{\frac{1}{q}}, \quad l_g^{q+}(\phi)(x)= \limsup_{n \to +\infty} \Big(  \prod_{k=n-q}^{n-1} \phi(g^k(x)) \Big )^{\frac{1}{q}}.$$

These maps verify the following inequalities:
$$0\leq l_g^-(\phi )(x)\leq l_g^+(\phi )(x)\leq +\infty,$$
$$0\leq L_g^-(\phi )(x)\leq L_g^+(\phi )(x)\leq +\infty, $$
$$ 0\leq l_g^{q-}(\phi )(x)\leq l_g^{q+}(\phi )(x)\leq +\infty$$

\medskip
\noindent for all $x\in X$. We will deal with $l_g(\phi)$,  $L_g(\phi )$, $l_g^q(\phi)$, whose domains are $\text{dom}(l_g(\phi ))=\{x \in X\hspace{2pt}|\hspace{2pt} l_g^-(\phi)(x)= l_g^+(\phi)(x)\}$ and similarly for $L_g(\phi)$ and $l_g^q(\phi)$.

\begin{definition}\label{Lyafunction} Let $(X,g)$ be a discrete semi-flow and $\phi \colon X \to [0,+\infty)$ a continuous map. The \emph{Lyapunov function} of $\phi$ is the function $L_g(\phi)$ described as
$$L_g(\phi)(x) =  \lim_{n\to +\infty} \Big( \prod_{k=0}^{n-1} \phi(g^k(x)) \Big )^{\frac{1}{n}}.$$
\end {definition}

\begin{remark} Note that the domain of $ L_g(\phi)$ is a subset of $X$ while the codomain is contained in $[0, +\infty]$.
Moreover, if $\phi \colon X \to [0, +\infty)$ is a bounded map, then the codomain of $ L_g(\phi)$ is a bounded subset of $[0, +\infty)$.
\end{remark}

\begin{remark} Lyapunov exponents have been used to study  local stability problems.
These are computed by taking time average functions (limits of arithmetic means) which are related with the convergence or divergence of nearby  orbits. 
In this paper in order to avoid some underflows,  which can appear when computing the limits of arithmetic means, and  using  the standard exponential function arithmetic means are transformed into   geometric means and we compute certain limits of these geometric means. That is,  instead of computing Lyapunov exponents we compute functions or constants  whose logarithms are Lyapunov exponents.

The sensitive dependence on initial conditions can be detected with the aid of Lyapunov functions or Lyapunov exponents. In our case, when the value of a  Lyapunov function at a point  is less than 1, we have local stability at the point and the nearby orbits converge to an attracting cycle. Notice that these Lyapunov functions are different from other type of Lyapunov functions which are used to detect local stability on continuous dynamics.

\end{remark}


\begin{definition} Let $(X,g)$ be a discrete semi-flow and $\phi \colon X \to [0,+\infty)$ a continuous map. Given $r, s \in \N$, $r\leq s$, the \emph{$[r,s]$-aproximation} of $L_g(\phi)$ is the map $L_g^{[r,s]}(\phi) \colon X \to \R_+$ described as
$$L_g^{[r,s]}(\phi) (x)=\left (\prod_{k=r}^s \phi (g^k(x))\right )^{\frac{1}{s+1-r}}. $$
Observe that the domain of $L_g^{[r,s]}(\phi)$ is $X$.
\end {definition}

\begin{definition} \label{fphicondition} Let $(X,g)$ be a discrete semi-flow and $\phi \colon X \to \R_+$ a continuous map.
Then  $(X,g, \phi)$ is said to verify the \emph{finite-nullity condition} if $g$ is a surjective open map with finite fibres, $\phi^{-1}(0)$ is a finite subset and $\phi$ is a bounded map.
If $\varphi \colon X\to \{-\infty\}\cup \R$ is a continuous map, then  $(X,g, \varphi)$ is said to verify the \emph{finite-infinite condition} if $g$ is a surjective open map with finite fibres, $\varphi^{-1}(-\infty)$ is a finite subset and $\varphi$ is an upper bounded map.
\end{definition}

\begin{definition} Let $(X,g)$ be a discrete semi-flow. The left $g$-invariant subset determined by $X'\subset X$ is
the least left $g$-invariant subset ${\rm{linv}}(X')$ of $X$ containing $X'$. It is given by
$${\rm{linv}}(X')=\bigcup_{k=0}^{\infty} g^{-k}(X').$$
\noindent If $\mathcal{P}(X)$ denotes the power set of $X$, then $g^{-1}\colon \mathcal{P}(X)\to \mathcal{P}(X)$ stands for the inverse image map, $g^{-k}= g^{-1} \circ\cdots^{(k)} \circ g^{-1} $ the $k$-th composition and $g^0=id_{\mathcal{P}(X)}$.
\end{definition}

\begin{proposition} \label{finitenullity}  Let $(X,g)$ be a discrete semi-flow and $\phi \colon X \to \R_+$ a continuous map.
If $(X,g, \phi)$ verifies the finite-nullity condition, then the following statements hold true:

\begin{itemize}
\item [(i)] If $x\in {\rm{linv}}(\phi^{-1}(0))$, then there exists $j_x\in \N$ such that $g^{j_x}(x) \in  \phi^{-1}(0)$ and $$L_g(\phi )(x)=\lim_{n\to +\infty} L_g^{[j,n]}(\phi ) (x)=0\hspace{10pt}\mbox{for}\hspace{4pt}j\leq j_x.$$

\item [(ii)] If $x \in {\rm{linv}}(\phi^{-1}(0))$ and there exists $j_x\in \N$ satisfying $g^{j_x}(x) \in  \phi^{-1}(0)$, $g^{j}(x) \not \in \phi^{-1}(0)$ for $j>j_x$ and $g^{j_x+1}(x) \in {\rm{dom}}(L_g(\phi))$, then $$L_g(\phi )(g^j(x))=\lim_{n\to +\infty} L_g^{[j,n]}(\phi )(x)\hspace{10pt}\mbox{for}\hspace{4pt}j>j_x.$$

\item [(iii)] If $x \not \in {\rm{linv}}(\phi^{-1}(0))$ and $x \in {\rm{dom}}(L_g(\phi ))$, then
$$L_g(\phi )(x)=\lim_{n\to +\infty} L_g^{[j,n]}(\phi )(x)\hspace{10pt}\mbox{for all}\hspace{4pt}j \in \N.$$
\end{itemize}
\end{proposition}

Correspondingly, one can develop a similar version using arithmetic means. Indeed, given any map $\varphi \colon X \to \R $ and bearing in mind that the infimum $\text{inf}\, T$ and the supremum $\text{sup}\, T$ of any nonempty subset $T\subset \R$, are elements in $\{-\infty\}\cup \R \cup\{+\infty\}$, we can consider the following functions:

 $$ a_g^- (\varphi) (x) =  \liminf_{n\to +\infty}\, (  \varphi\left(g^n (x)\right) ), \quad a_g^+ (\varphi) (x) = \limsup_{n\to +\infty}\, ( \varphi\left(g^n (x)\right) ),$$

 $$ A_g^- (\varphi) (x) =  \liminf_{n\to +\infty}\,( \frac{1}{n} \sum_{k=0}^{n-1} \varphi\left(g^k (x)\right)), \quad A_g^+ (\varphi) (x) = \limsup_{n\to +\infty}\, ( \frac{1}{n} \sum_{k=0}^{n-1} \varphi\left(g^k(x)\right)), $$

  $$ a_g^{q-} (\varphi) (x) = \liminf_{n\to +\infty}\, ( \frac{1}{q} \sum_{k=n-q}^{n-1} \varphi\left(g^k (x)\right) ),  \quad a_g^{q+} (\varphi) (x) = \limsup_{n\to +\infty}\, ( \frac{1}{q} \sum_{k=n-q}^{n-1} \varphi\left(g^k (x)\right) ).$$

We can take $\text{dom}\, a_g( \varphi)=\{x \in X \hspace{2pt}|\hspace{2pt} a_g^-(\varphi)(x)=a_g^+(\varphi)(x)\}$ and similarly the corresponding domains $\text{dom}\, A_g \varphi)$, $\text{dom}\, a_g^q( \varphi)$. Then we obtain the functions $a_g( \varphi), A_g( \varphi)$ and $a_g^q( \varphi)$ defined on their corresponding domains. Notice that the definition of $a_g(\varphi)$ is similar to $l_g(\phi)$, the main difference being that the codomain of $a_g(\varphi)$ is $\{-\infty\}\cup \R \cup\{+\infty\}$ instead of $[0,+\infty]$.

\begin{remark} Note that the domain of $ A_g(\varphi)$ is a subset of $X$ and the codomain is contained in $\{-\infty\}\cup \R \cup\{+\infty\}$. Also observe that if $\varphi \colon X \to \mathbb{R}$ is upper bounded, then the codomain of $ A_g(\varphi)$ is an upper bounded subset of $\{-\infty\}\cup \R $.
\end{remark}

Now we recall the usual notion of Time Average:

 \begin{definition}Let $(X,g)$ be a discrete semi-flow and $\varphi \colon X \to \R$ a continuous map. The \emph{time average} of $\varphi$ is the function $A_g(\varphi)$ described as
$$A_g(\varphi)(x)=\lim_{n\rightarrow +\infty}\;   \frac{1}{n} \sum_{k=0}^{n-1} \varphi\left(g^k (x)\right).$$
\end {definition}

On measurable spaces the arithmetic means can be extended to the space average:

\begin{definition} Let  $(X,\Sigma)$ be a measurable space and $\mu$ a measure where $\mu(X)$ is finite and nonzero. The \emph{space average} $S(\varphi)$  of an integrable function $\varphi\colon X \to \R$ is described as
$$S(\varphi) =\frac 1{\mu(X)} \int_X \varphi\,d\mu .$$
\end{definition}

In our study of attracting basins it will be crucial the use of discrete measures. Consider a measurable space $(X,\Sigma )$ and two measures $\mu$, $\nu$ on it. Then $\mu$ is said to be \emph{discrete with respect to} $\nu$ if there exists a countable subset $S$ of $X$ such that for every $s\in S$, the singleton $\{s\} \in \Sigma$, $\nu (S)=0$ and  $\mu (X\backslash S)=0.$

Now, given a countable subset $S=\{s_0, s_1, s_1,\cdots\}\subset X$ such that  $\{s\} \in \Sigma$ for every $s\in S$, one can easily check  that $\mu$ is  discrete with respect to $\nu$ on  $S$ if, and only if, there is a sequence $a_{0}, a_{1},a_{2},\dots$ of numbers in $[0,+\infty ]$ such that $\mu =\sum _{i}a_{i}\delta _{s_{i}}$ and $\nu(S)=0$. Here $\delta _{s_{i}}$ denotes the so-called \emph{Dirac measure}, defined as
$$
\delta _{s_{i}}(A)=
\begin{cases} 1, &  {\mbox{ if }}s_{i}\in A, \\
0,  & {\mbox{ if }}s_{i}\not \in A. \end{cases}
$$

\begin{remark} Given a measure space $(X,\Sigma, \mu )$ with $0\not = \mu (X) <+\infty$, a continuous map $g \colon X \to X$ and an integrable function $\varphi \colon X \to \R$, the time average and space average may be different. However, if the map $g$ is ergodic (i.e.,  the conditions $B\in \Sigma$ and $g^{-1}(B)=B$ imply $\mu(B)=\mu(X)$ or $\mu(B)=0$) and the measure $\mu$  is invariant by $g$ (i.e., $\mu(g^{-1}(B))=\mu(B), B \in \Sigma$), then  by  the named   Birkhoff Ergodic Theorem \cite{Birkhoff},  the time average is equal to the space average  $A_g(\varphi)=S(\varphi)$ almost everywhere.  For a version of the notion of  ergodic linear operators we refer the reader to \cite{grosse2011linear}.

\end{remark}

\begin{definition}Given a measure space $(X,\Sigma, \mu )$ with $0\not = \mu (X) <+\infty$ and  a continuous  map $g \colon X \to X$, we say that $(X,\Sigma, \mu, g )$ is a \emph{Birkhoff space}, or just a \emph{$B$-space}, if for any continuous function $\varphi \colon X \to \R$ the following equality
 \begin{equation}
\label{Birkoffformula}
\lim_{n\rightarrow+\infty}\; \frac{1}{n} \sum_{k=0}^{n-1} \varphi \left(g^k(x)\right) = \frac 1{\mu(X)}\int_X \varphi\,d\mu
\end{equation}
\noindent holds almost everywhere. In case of equality for every $x\in X$ we say that $(X,\Sigma, \mu, g )$ is a \emph{full $B$-space}.(If we consider integrable maps $\varphi \colon X \to \R$ instead of continuous maps, then we obtain the notion of \emph{$B^1$-space}). In some cases, we also have a different measure $\nu$ defined on $\Sigma$ and the formula \eqref{Birkoffformula} is satisfied up to a subset whose $\nu$-measure is zero; in this case we say that $(X,\Sigma, \mu,\nu, g )$ is a \emph{$\nu$-$B$-space} (or a \emph{$\nu$-$B^1$-space}).
\end{definition}

\begin{remark}\label{quasi} Given two measure spaces $(X_1,\Sigma_1, \mu_1 )$, $(X_2,\Sigma_2, \mu_2 )$, continuous  maps $g_1 \colon X_1 \to X_1$,
$g_2 \colon X_2 \to X_2$,  and continuous   (integrable)  functions $\varphi _1\colon X_1 \to \R$, $\varphi _2\colon X_2 \to \R$, then one can consider the sum $X=X_1+ X_2$, $g=g_1+ g_2$ and $\varphi=\varphi_1+ \varphi_2$ with the $\sigma$-algebra  $\Sigma$ generated by $\Sigma_1$ and $\Sigma_2$. In this way we obtain a new measure space $(X,\Sigma, \mu )$ and the sum measure $\mu=\mu_1+\mu_2$. Similarly, for a finite number of measure spaces $(X_i,\Sigma_i, \mu_i )$, $i\in \{1,...,m\}$,  we have the corresponding sum measure space $X=\sum_{i=1}^mX_i$. Notice that if each measure $\mu_i$ is invariant under $g_i$, then $\mu= \sum_{i=1}^m\mu_i$ is invariant under $g=\sum_{i=1}^m g_i$.
If a measure discrete semi-flow $(X,\Sigma, \mu, g )$ can be  expressed as a finite sum of  B-spaces ($B^1$-spaces)   $(X_i,\Sigma_i, \mu_i, g_i )$, then, for each $i \in \{1, \cdots, m\}$ and for every continuous (integrable) function $\varphi :X\rightarrow \mathbb{R}$ we have that
$$\lim_{n\rightarrow +\infty}\; \frac{1}{n} \sum_{k=0}^{n-1} \varphi \left(g^k (x)\right) = \frac 1{\mu(X_i)}\int_{X_i}\varphi\,d\mu $$
\noindent almost everywhere in $X_i$. (In some cases, we shall work with two measures $\mu_i, \nu_i$ and we can consider similar finite sums of $\nu_i$-$B$-spaces.)
\end{remark}

\subsection{The complex projective line and the Hopf fibration} \label{structCinf}

With the aim of studying the iteration of a complex rational map, we will use the following models of the Riemann sphere:
$S^2=\{(r_1,r_2,r_3) \in \R^3 \mid r_1^2+r_2^2+r_3^2=1\}$ the unit 2-sphere,
$\hat \C = \C \cup \{\infty\}$ the Alexandroff compactification of $\C$ and
$\p^{1}(\C)$ the complex projective line, which is given by
\begin{equation}
\label{projective}
\p^{1}(\C)=\big ( \C^2\setminus \{(0,0)\}\big )/ \sim,
\end{equation}
where $(z,t)\sim(z',t')$ if, and only if, there exists $\lambda\in\C\setminus\{0\}$ such that $(z,t)=(\lambda z',\lambda t')$
(the equivalence class of $(z,t)$ is denoted by $[z\colon t]$).

There are homeomorphisms $\tilde\theta\colon S^2\to\hat \C$ and $\theta\colon \p^1(\C)\to\hat \C$, given as
$$\tilde\theta(r_1,r_2,r_3)=
\begin{cases}
\frac{r_1}{1-r_3}+i\frac{r_2}{1-r_3}, & \text{if }r_3<1, \\
\infty, & \text{if }r_3=1.
\end{cases}
 \quad ;
 \quad \quad
\theta([z\colon t])=\begin{cases}
z/t, & \text{if }t\ne0, \\
\infty, & \text{if }t=0.
\end{cases}
$$
\noindent and the induced composite homeomorphism, $\tilde \theta^{-1}\circ \theta\colon \p^1(\C)\to S^2$, given as
\begin{equation}
\label{isomorphism}
(\tilde\theta^{-1}\circ \theta)([z:t])) =
\left(\frac{\bar zt+z\bar t}{\bar tt+z\bar z},\frac{i(\bar zt-z\bar t)}{\bar tt+z\bar z},\frac{-\bar tt+z\bar z}{\bar tt+z\bar z}\right)
\end{equation}
%
where $\bar z$ denotes the conjugate of $z$. Note that the chordal (euclidean) metric on $S^2$ and the homeomorphism in \eqref{isomorphism} induce a \emph{chordal metric} on $\p^{1}(\C)$ and on $\hat \C.$ Therefore we are actually dealing with isometries; moreover, $\theta$ and $\tilde \theta$ are isomorphisms in the category of Riemann surfaces and complex analytic maps.

\bigskip

In this paper we also contemplate the 3-sphere $S^3$ as a subspace of $\C^2$:
\begin{equation}
\label{3sphere}
S^3=\{(z,t) \in \C^2 \hspace{3pt}|\hspace{3pt} \, \abs{z}+\abs{t}=1\}
\end{equation}
where $|\lambda|$ denotes the absolute value (or modulus) of a complex number $\lambda$. We are working with a composite of quotients:
$$\xymatrix{\C^2 \setminus \{(0,0)\} \ar[r]^(.65)p & S^{3} \ar[r]^(.4)q & \p^{1}(\C)}.$$
Here, the surjective map $p \colon \C^2 \setminus \{(0,0)\} \to S^{3}$ is given as $p(z,t)=[z ,t], $ with
\begin{equation}
\label{Normalization}
[z,t]=(\frac{z}{|z|+|t|}, \frac{t}{|z|+|t|}), \quad (z,t)\in \C^2 \setminus \{(0,0)\},
\end{equation}

\noindent whereas the surjective map $q \colon S^{3} \to \p^{1}(\C)$ is given as
\begin{equation}
\label{esetrestoplane}
 \quad q(z,t)=[z :t], \quad (z,t) \in S^{3}.
\end{equation}

These surjections induce the following equivalence relations:

\medskip
$\bullet $ Let $(z,t), (z',t')$ be in $\C^2 \setminus \{(0,0)\}$. Then $(z,t)\sim_p (z',t')$ if, and only if, $(z',t')= \lambda (z,t)$ for some real $\lambda >0$.

$\bullet $ Let $(z,t), (z',t')$ be in $S^{3}$. Then $(z,t)\sim_q (z',t')$ if, and only if, $(z',t')= \lambda (z,t)$ for some $\lambda \in \C, \abs{\lambda}=1$.

\medskip

\begin{definition}\label{normalizedhc}
Any representative $(z,t)$ of a point $[z\colon t]\in \p^{1}(\C)$ is called \emph{homogeneous coordinates} of this point. In particular we shall call the \emph{normalized homogeneous coordinates} to the pair
$[z,t]=\big (\frac{z}{|z|+|t|},  \frac{t}{|z|+|t|}\big )$ (see formula  \eqref{Normalization}).
\end{definition}

\begin{definition} The map  $q\colon S^{3}\to \p^{1}(\C)$ is said to be the \emph{Hopf fibration}, and
the composite $q\circ p\colon \C^2 \setminus \{(0,0)\}  \to  \p^{1}(\C)$ the \emph{expanded-Hopf fibration}.
\end{definition}

\begin{remark}
The Hopf fibration \cite{Hopf} was applied to prove that the 3rd homotopy group of the 2-sphere is not trivial. In this paper,   Hopf fibration is used to avoid overflow and underflows problems during the iteration of a rational map.
\end{remark}

\subsection{Complex rational maps} \label{rf}

\subsubsection{Homogeneous pairs and  complex rational maps}

Let $\C[x]$ be  the ring of polynomials with complex coefficients and $\C(x)$ the field of complex rational functions
$f(x)=\frac{A(x)}{B(x)}$, where $A, B \in \C[x]$, $B$ is not the  zero polynomial and the equations $A(x)=0$, $B(x)=0$ do not have a common root solution.
The degree $d$ of $f$ is given by $d=\max\{\text{degree}(A), \text{degree}(B)\}$; the degree of the zero polynomial is taken as $-\infty$.

\medskip
Let $\C[z,t]$ be the ring of bivariate polynomials with complex coefficients. Note that if we have a quotient  $\frac{A(x)}{B(x)}$ and $d=\max\{\text{degree}(A), \text{degree}(B)\}$, taking $x=z/t$, we can consider the bivariate polynomials $F(z,t), G(z,t) \in \C[z,t]$ given as
\begin{equation}
\label{rationaltohomogeneous}
F(z,t)=t^d A(z/t), \quad G(z,t)=t^d A(z/t).
\end{equation}

\begin{definition}\label{bivariate} A bivariate polynomial $F \in \C[z,t]$ is said to be \emph{homogeneous}
if $F=0$ or there exists an integer $k\geq 0$ such that $F(\lambda z, \lambda t)=\lambda^k F(z,t)$ for every $z, t \in \C$ and $\lambda \in \C\setminus \{0\}$. One defines $\text{degree} (F)=k$ for this latter case, whereas $\text{degree} (F)=-\infty$ for the case $F=0$.
\end{definition}

We shall denote as $\C_h[z,t]$ the subset of $\C[z,t]$ consisting of all homogeneous bivariate polynomials. Any homogeneous bivariate polynomial $F$ with degree $k\geq 0$ can be expressed as $F(z,t)=\sum_{i=0}^k a_i z^i t^{k-i}$ where each $a_i \in \C$. However, if
$n=\min\{m \hspace{3pt}|\hspace{3pt} a_i=0 \text{ for } i> m\}$, then we will write $F$ as
$$F(z,t)=\sum_{i=0}^n a_i z^i t^{k-i}$$ \noindent where each $a_i \in \C,$  $a_n \not = 0$ and $n\leq k.$

\medskip

Now if we consider $\C[x]  \underline{\times} \N:=\{(A,d)\in  \C[x]  {\times} \N\hspace{3pt}|\hspace{3pt} \text{degree}(A)\leq d\},$
then we can  take the \emph{\it homogenization operator} ${\rm H}: \C[x] \underline{\times} \N \to \C_h[z,t]$, defined as follows:
\begin{equation}
\label{homogeneousversion}
{\rm H}(A,d)=\begin{cases}
0, & \text{if } A=0, \\
\sum _{i=0}^na_iz^it^{d-i}, & \text{if }
A=\sum _{i=0}^na_ix^i\neq 0, a_n\not =0.
\end{cases}
\end{equation}

We note that for $A\not =0$ we have ${\rm H}(A,d)(\lambda z, \lambda t)=\lambda^d {\rm H}(A,d)(z, t)$.

\begin{definition}\label{homogeneouspair}  Given
$F,G \in\C_h[z,t]$, the pair $(F,G)$ is said to be a \emph{homogeneous pair}, or just \emph{homogeneous}, if either $FG=0$ (i.e., their product vanishes) or $F$ and $G$ have the same degree. In any case,
$$d=\max\{\text{degree}(F),\text{degree}(G)\}$$ \noindent is called \emph{the degree of} $(F,G).$
When $(F,G)$ is homogenous with $G\not =0$, then $(F,G)$ is said to be \emph{r-homogeneous}.
\end{definition}

We will consider the following notation:
\begin{equation}
\label{homogeneouspairs}
\C_h[z,t]\underline{\times} \C_h[z,t]:= \{(F,G) \in \C_h[z,t]{\times} \C_h[z,t]\hspace{2pt}|\hspace{2pt} (F,G) \textit {\rm{ is homogeneous}}\},
\end{equation}
\vspace{-15pt}
\begin{equation}
\label{rhomogeneouspairs}
\C_h[z,t]\underline{\times}_{\rm{r}} \C_h[z,t]:= \{(F,G) \in \C_h[z,t]{\times} \C_h[z,t]\hspace{2pt}|\hspace{2pt} (F,G) \textit {\rm{ is r-homogeneous}}\}
\end{equation}

It is also interesting to analyze the case of $\frac{A(x)}{B(x)}$ being irreducible. In this case the corresponding homogeneous pair $(F,G)$ is an  r-homogeneous irreducible pair:

\begin{definition}\label{irreducible} If $(F,G)\in \C_h[z,t]\underline{\times} \C_h[z,t] $, then $(F,G)$ is said to be \emph{irreducible} if either $(F,G)=(0,0)$ or $(F,G)\not =(0,0)$ and satisfies the following property: $$(F,G)=(HF_1,HG_1)\Rightarrow \text{degree}(H)=0.$$ The pair $(F,G)$ is said to be \emph{r-irreducible} if  $(F,G)$ is an r-homogeneous irreducible pair.
\end{definition}

Accordingly, we consider the following notation:
\begin{equation}
\label{ihomogeneouspairs}
\C_h[z,t]\underline{\times}_{\rm{i}} \C_h[z,t]:= \{(F,G) \in \C_h[z,t]\underline{\times} \C_h[z,t]\hspace{2pt} |\hspace{2pt} (F,G) \textit {\rm{ is irreducible}}\},
\end{equation}
\vspace{-15pt}
\begin{equation}
\label{rhomogeneouspairs}
\C_h[z,t]\underline{\times}_{\rm{ri}} \C_h[z,t]:= \{(F,G) \in\C_h[z,t]\underline{\times} \C_h[z,t]\hspace{2pt}|\hspace{2pt} (F,G) \textit {\rm{ is r-irreducible}}\}.
\end{equation}

\medskip
Now, given any pair $R=(F,G)\in C_h[z,t]\underline{\times}_{\rm{ri}} \C_h[z,t]$ we get an induced commutative diagram:
$$\xymatrix{ \ar[d]_p \C^2 \setminus \{ (0,0) \} \ar[r]^R & \C^2\setminus \{ (0,0)\} \ar[d]^p \\
\ar[d]_q S^{3} \ar[r]^{R^S}  & S^{3} \ar[d]^q \\
\p^{1}(\C) \ar[r]^{R^P} &\p^{1}(\C) }
$$
where $R$, $R^S$ and $R^P$ are defined as follows:
\begin{equation}
\label{erre}
R(z,t):=(F(z,t),G(z,t)),
\end{equation}
\vspace{-20pt}
\begin{equation}
\label{ese}
R^S(z,t):=[F(z,t),G(z,t)],
\end{equation}
\vspace{-20pt}
\begin{equation}
\label{pe}
R^P([z:t]):=[F(z,t):G(z,t)]
\end{equation}
where $[F(z,t),G(z,t)]=(\frac{F(z,t)}{|F(z,t)|+|G(z,t)|}, \frac{G(z,t)}{|F(z,t)|+|G(z,t)|})$, see formula \eqref{Normalization}.

\begin{definition} Given a homogeneous pair $R=(F,G) \in \C_h[z,t]\underline{\times}_{\rm{ri}} \C_h[z,t]$, the pair
$(R^S, R^P)$ is said to be  the \emph{Hopf-endomorphism} induced by $R$  and the triple
$(R, R^S, R^P)$ is the \emph{expanded-Hopf-endomorphism} induced by $R$.
\end{definition}

\subsubsection{Normalization of a homogeneous pair}\label{normphi}

In the following paragraphs we are considering what we call the normalization of homogeneous pairs. The use of normalized  homogeneous pairs  will avoid computational overflows in our algorithms.


\begin{definition}\label{normapar} Let $F\in \C_h[z,t]$ be a homogeneous polynomial  given by $F(z,t)=a_0 t^d+a_1 z t^{d-1}+\cdots +a_n z^n t^{d-n}$.  We define   the \emph{norm} $\norm{F}$ as
$$\norm{F} = |a_0|+|a_1|+\cdots +|a_n|. $$
The \emph{norm} of $(F,G)\in \C_h[z,t]\underline{\times} \C_h[z,t]$ is defined as $$\norm{(F,G)} :=\norm{F}+\norm{G}.$$ Finally, the \emph{normalization} $[F,G]$ of a pair $(F,G)\in \C_h[z,t]\underline{\times} \C_h[z,t]$ is given as follows:
 \begin{equation}
\label{Normalizationpair}
 [F,G]:=\begin{cases}
(0,0), & \text{if } (F,G)=(0,0), \\
(\frac{1}{\norm{(F,G)} }F,\frac{1}{\norm{(F,G)}} G), & \text{if } (F,G)\neq (0,0).
\end{cases}
\end{equation}

\noindent A homogeneous pair $(F,G)$ is said to be \emph{normalized} if $[F,G]=(F,G)$.
 \end{definition}

Denote as $\big (\C_h[z,t]\underline{\times} \C_h[z,t]\big )_{\text{n}}$ the subset of $\big (\C_h[z,t]\underline{\times} \C_h[z,t]\big )$ consisting of all normalized homogeneous pairs. Consider also the notation
$$\begin{array}{l}
D^2= \{z\in \C\hspace{2pt} |\hspace{2pt} \abs{z} \leq 1\}, \\
D^4=\{(z,t) \in \C^2\hspace{2pt}|\hspace{2pt} \abs{z}+\abs{t}\leq 1\}.
\end{array}$$
If $(a_0,a_1, \cdots, a_n)$ and $(b_0,b_1, \cdots, b_m)$ are the respective coefficients of a homogeneous pair $(F,G)$ of degree $d$  and
$(z,t) \in D^2 \times D^2$, then one has
\begin{center}$|a_0 t^d+a_1 z t^{d-1}+\cdots +a_n z^n t^{d-n}| \leq  |a_0|+|a_1|+\cdots +|a_n|,$\end{center}
\begin{equation}\label{bounded}
|b_0 t^d+b_1 z t^{d-1}+\cdots +b_m z^m t^{d-m}|\leq  |b_0 |+| b_1|+\cdots +|b_m|.
\end{equation}

\begin{proposition}\label{disks} If $R=(F,G) \in \C_h[z,t]\underline{\times} \C_h[z,t]$ is a normalized homogeneous pair, then we have a chain of inclusions:
$$R(S^{3})\subset R(D^{4})\subset  R\big (  D^2 \times D^2 \big )\subset D^4 \subset   D^2 \times D^2.$$
Furthermore, if $R=(F,G)$ is a normalized r-irreducible pair, then we obtain the following commutative diagram:
$$\xymatrix{ D^4 \setminus \{(0,0)\} \ar[rr]^{R=(F,G)} & & D^4 \setminus \{(0,0)\} \ar[d]_{p} & \\
\ar@{^{(}->}[u]^j  \ar[urr]^{R\circ j}\ar[d]_q S^{3} \ar[rr]_{R^S=[F,G]} &   & \ar[d]_q S^{3}\\ \
\p^{1}(\C) \ar[rr]_{R^P=[F:G]} & & \p^{1}(\C)   &}
$$
where $j\colon S^{3}\hookrightarrow D^4\setminus \{(0,0)\}$ is the canonical inclusion.
 \end{proposition}

\subsection{Computing the spherical derivative of a rational function}\label{computingsphericalderivative}

Now we develop a method that permits an easy implementation of the spherical derivative avoiding the problems associated to poles and the definition at the infinite point. Such problems can be found in the following standard definition of spherical derivative:

 \begin{definition}\label{sphericalderivative}
Given a rational map $f \colon \hat \C \to \hat \C$, the \emph{spherical derivative} of $f$ is the map
$f^\# \colon \hat \C \to \R_+$ defined as
$$ f^\# (x_0) = \abs{f'(x_0) } \frac{1+\abs{x_0}^2}{1+\abs{f(x_0)}^2}$$
where $f'(x_0)$ is the derivative of $f$ at $x_0$. If the right expression is not defined at $x_0$, then $f^\# (x_0) $ is taken as the limit
$$\lim _{x\rightarrow x_0}\abs{f'(x) } \frac{1+\abs{x}^2}{1+\abs{f(x)}^2}$$ A point $x_0 \in \hat \C$ it is said to be \emph{critical} when $f^\# (x_0)=0$.
\end{definition}

Suppose a rational map $f(x)=A(x)/B(x)$ and consider its associated homogeneous pair $(F,G)$:
$$
\begin{array}{l}
F(z,t)=t^d A(z/t), \\
G(z,t)=t^d   B(z/t)
\end{array}$$ \noindent where $d=\max\{\text{degree}(A),\text{degree}(B)\}$. The standard derivative of the map $f \colon \hat \C \to  \hat  \C$ is given as $$f'(x)=\frac{A'(x)B(x) -A(x) B'(x)}{B(x)^2}.$$
Now, given any polynomial $A(x)=\sum_{k=0}^n a_k x^k$, the conjugate of $A(x)$ is plainly given as $\overline{A(x)}= \sum_{k=0}^n {\bar a}_k {\bar x}^k$. If $\bar A$ stands for the polynomial whose coefficients are $(\bar a_0, \cdots, \bar a_n)$, then one has $\overline{A(x)}={\bar A}(\bar x)$.

Paying regard to the considerations above we have for $f(x)=A(x)/B(x)$ and $x=z/t$ that:

\begin{eqnarray}
f'(x) & = & \frac{A'(z/t) B(z/t) -A(z/t) B'(z/t)}{B(z/t)^2}, \\
\frac{1+\abs{x}^2}{1+\abs{f(x)}^2} & = &  \frac{(1+ (z/t) (\bar z/\bar t))\hspace{1pt} B(z/t) \bar B(\bar z/\bar t) }{B(z/t) \bar B(\bar z/\bar t)+A(z/t) \bar A(\bar z/\bar t)}.
\end{eqnarray}

Consequently,
$$f'(x) \frac{1+\abs{x}^2}{1+\abs{f(x)}^2} =\frac{(A'(z/t) B(z/t) -A(z/t) B'(z/t))\hspace{1pt}(1+ (z/t) (\bar z /\bar t))}{B(z/t) \bar B(\bar z/\bar t)+A(z/t) \bar A(\bar z/\bar t)} \frac{\bar B( \bar z/ \bar t) }{B(z/t)}.$$


%
But, clearly, $\abs{\frac{\bar B( \bar z/ \bar t) }{B(z/t)}}=1$ and therefore
\begin{equation}
\label{finaleq}
f^\# (x) = \left |\frac{(A'(z/t) B(z/t) -A(z/t) B'(z/t)) (1+ (z/t) (\bar z /\bar t))}{B(z/t) \bar B(\bar z/\bar t)+A(z/t) \bar A(\bar z/\bar t)}\right |.
\end{equation}


\medskip
Take $d=\max\{\text{degree}(A),\text{degree}(B)\}$, $e:=\text{degree}(A)+ \text{degree}(B)$ and \linebreak $e^*:=\text{degree}(A'B-AB')$. Then we have the following interesting property whose proof consists only of doing a simple check:

\begin{lemma}\label{lemainterest} Suppose $A(x)=\sum_{i=0}^n a_i x^i$ with $a_n \not =0$ and $B(x)=\sum_{j=0}^m b_j x^j$ with $b_m\not = 0$ (so that $\text{degree}(A)=n$ and $\text{degree}(B)=m$). The following holds true:
\begin{itemize}
\item[(i)] If $n\not =m$, then $e^*=e-1$ and the leading coefficient of $A'B-AB'$ is equal to $(n-m)a_n b_m$.
\item[(ii)] If $n=m$, then $e^*<e-1$.
\end{itemize}
In both cases, $2d > e^*+1$.
%
\end{lemma}

\bigskip
Let $A,B$ be polynomials such that $A B \not = 0$ and remember that we are denoting $d=\max\{ \text{degree} (A), \text{degree} (B) \}$.
Now we consider the maps $u, v \colon \C \times (\C \setminus \{0\} )\to \C$ given as
\begin{equation}
\label{definitionu}
u(z, t):= t^d {\bar t}^d  (A'(z/t) B(z/t) -A(z/t) B'(z/t)) (1+ (z/t) (\bar z /\bar t)),
\end{equation}
\vspace{-20pt}
\begin{equation}
\label{definitionv}
v(z,t):= t^d {\bar t}^d ({B(z/t) \bar B(\bar z/\bar t)+A(z/t) \bar A(\bar z/\bar t))}.
\end{equation}

We have the following straightforward equations:
\begin{equation}
\label{descompuesto}
 t^d {\bar t}^d= \abs{t}^{2d-(e^*+2)} \abs{t}^{e^*} t \bar{t},
\end{equation}
\vspace{-25pt}
\begin{equation}
\label{derivadas}
A'(x)=\sum_{i=0}^{n-1} (i+1) a_{i+1} x^{i}, \quad B'(x)=\sum_{j=0}^{m-1} (j+1) b_{j+1} x^{j}.
\end{equation}


Let $\text{coef}_k$ denote the $k$-th coefficient of the polynomial $A'B-AB'$. As a consequence of equation \eqref{derivadas} and the definition of $e^*$,  we get the following equations, which hold for every $z\in \C$ and every $t \in \C\setminus \{0\}$:
\begin{equation}
\label{emenosunoa}
 t^{e^*} (A'(z/t) B(z/t)-A(z/t) B'(z/t))= \sum_{k=0}^{e^*} \text{coef}_k\hspace{1pt}z^{k} t^{e^*-k},
\end{equation}
\vspace{-16pt}
\begin{equation}
\label{unoconj}
 t  {\bar t} (1+ (z/t) (\bar z /\bar t))= t \bar t + z \bar z,
\end{equation}
\vspace{-16pt}
\begin{equation}
\label{BBconj}
 t^{d}{\bar t}^{d}\hspace{1pt} B(z/t) {\bar B}(\bar z/ \bar t)= \sum_{i=0}^{m} \sum_{j=0}^{m} b_{i} {\bar b}_j z^{i} t^{d-i} {\bar z}^{j} {\bar t}^{d-j},
\end{equation}
\vspace{-16pt}
\begin{equation}
\label{AAconj}
 t^{d}  {\bar t}^{d}\hspace{1pt}(A(z/t) {\bar A}(\bar z/ \bar t)= \sum_{i=0}^{n}  \sum_{j=0}^{n}a_{i} {\bar a}_j z^{i} t^{d-i} {\bar z}^{j} {\bar t}^{d-j}.
\end{equation}
Here, the right-hand side of equations \eqref{emenosunoa}, \eqref{unoconj}, \eqref{BBconj} \eqref{AAconj} represent continuous functions.

\begin{definition}\label{sphericalpair}
The \emph{spherical homogeneous pair $(U,V)$ induced by $A, B$} is given by
$$
\begin{array}{l}
U(z,t):= \abs{t}^{2d-(e^*+2)}  \left (\sum_{k=0}^{e^*} \text{coef}_k  z^{k} t^{e^*-k}\right ) \left ( t \bar t + z \bar z  \right ), \\
V(z,t):=\left (  \sum_{i=0}^{m} \sum_{j=0}^{m} b_{i} {\bar b}_j z^{i} t^{d-i} {\bar z}^{j} {\bar t}^{d-j}\right ) + \left (  \sum_{i=0}^{n}  \sum_{j=0}^{n}a_{i} {\bar a}_j z^{i} t^{d-i} {\bar z}^{j} {\bar t}^{d-j}
\right ).
\end{array} $$
\end{definition}

\begin{remark}\label{continuous} The pair $(U, V)$ represents a continuous function from $\C^2$ to $\C^2$.
Note that the family of all polynomials on variables $z$, $\bar z$, $t$, $\bar t$ is just an extension of the family of polynomials on $z, t$ together with the family of polynomials on $\bar z, \bar t$.
\end{remark}

The spherical derivative and its spherical homogeneous pair are related as follows. Remember that for a rational map $f(x)=\frac{A(x)}{B(x)}$ we have the notation $d=\max\{\text{degree}(A), \text{degree}(B)\}$, $e=\text{degree}(A)+\text{degree}(B)$ and $e^*=\text{degree}(A'B-AB')$.

\begin{proposition}\label{computinspherical} Let $f\colon \hat\C \to \hat\C$ be a rational map and take an irreducible representation $f(x)=\frac{A(x)}{B(x)}$. If we consider the maps $U, V \colon \C^2 \to \C$ given by Definition \ref{sphericalpair}, then the following statements hold true:
\begin{itemize}
\item[(i)]
$U(\lambda z,  \lambda t)=\lambda^d  {\bar \lambda}^d  U(z, t)$ and
$V(\lambda z,  \lambda t)=  \lambda^d  {\bar \lambda}^d  V(z, t)$, for all $\lambda \in \C$.

In particular, if $|\lambda |=1$, then $U(\lambda z,  \lambda t)=U(z, t)$ and $V(\lambda z,  \lambda t)=V(z, t)$.

\item[(ii)] $V(z, t) \not = 0$, for all $(z,t) \in \C^2\setminus \{(0,0)\}$.
     \item [(iii)] The induced maps
  $$S^3 \to \R ,\hspace{8pt} (z,t) \mapsto \abs{U(z,t)/ V(z,t)},$$
 $$\p^1(\C) \to \R, \hspace{8pt} [z:t]  \mapsto \abs{U(z,t)/ V(z,t)}$$ are   continuous.

 \item [(iv)] If $x=z/t$ for $(z,t) \in \C\times \C \setminus \{0\}$, then
 $$ f^\# (x) =\left | U(z,t)/ V(z,t) \right |, \quad  f^\# (\infty) = \left | U(z,0)/ V(z,0) \right |.$$

\end{itemize}
\end{proposition}

\begin{proof} ${}$

(i) is a consequence of Definition \ref{sphericalpair} and equation \eqref{unoconj}.

(ii) If $t\not =0$, since $\frac{A(x)}{B(x)}$ is an irreducible representation, we have $$V(z,t)= t^d {\bar t}^d ({B(z/t) \bar B(\bar z/\bar t)+A(z/t) \bar A(\bar z/\bar t))}\not = 0.$$

On the other hand, if $t=0$, then using the Definition \ref{sphericalpair}, we have
$$
V(z,0)=
\begin{cases}
 a_n \bar{a}_n   z^n {\bar z}^n,  & \text{ if } n >m, \\
  b_m \bar{b}_m  z^m {\bar z}^m, & \text{ if } m >n, \\
 (a_n \bar{a}_n +b_m \bar{b}_m)z^n {\bar z}^m,   & \text{ if } m=n.
\end{cases}
$$
This implies that $V(z,t)\not =0$, for every $(z,t) \in \C^2\setminus \{(0,0)\}$.

(iii) Taking into consideration (i) and (ii) above, one has that the induced maps $S^3 \to \R,$ $(z,t) \mapsto \abs{U(z,t)/ V(z,t)}$ and $\p^1(\C) \to \R,$   $[z:t]\mapsto \abs{U(z,t)/ V(z,t)}$ are well defined. As $p \colon \C^2 \setminus \{(0,0)\}\to S^{3}$, $q\colon S^{3}\to \p^{1}(\C)$ are quotient maps, the
continuity follows from Remark \ref{continuous}.

(iv) If $x=z/t \in \C$ with $t\not =0$, then, using equations \eqref{finaleq}-\eqref{AAconj} we have that
$$\left | U(z,t)/ V(z,t) \right |= \left | u(z,t)/ v(z,t) \right |= f^\# (x)
.$$

As a consequence of Definition \ref{sphericalderivative} and (iii) we also have that $f^\# (\infty) =\left |{U(z,0)}/{V(z,0)}\right |$.
\end{proof}

In the following corollary, we compute the value of the spherical derivative at the point at infinity. This value depends on the values of  $e^*$, $m$ and $n$ and $d=\max\{m, n\}$. See notation before  Lemma \ref{lemainterest}.

\begin{corollary} Using the notation above, we have
\begin{itemize}
  \item [(i)] If $2d>e^*+2$, then $U(z,0)=0$ and $f^\#(\infty)=0.$
  \item [(ii)] If $2d=e^*+2$, then $U(z,0)=\text{\rm coef}_{ e^*} z^{e^*}z {\bar z}$, that is,
$$
U(z,0)=
\begin{cases}
a_n b_m z^{n+m-1} z {\bar z},  & \text{ if } n =m+1,\\
 -a_n b_m  z^{n+m-1} z {\bar z},  & \text{ if } m =n+1, \\
\text{\rm coef}_{ e^*} z^{e^*} z {\bar z},  & \text{if } n=m.
\end{cases}
$$
In this case:
    $$
f^\#(\infty)=\abs{U(z,0)/V(z.0)}=
\begin{cases}
\abs{b_m}/\abs{{a}_n}, & \text{ if } n=m+1, \\
\abs{a_n}/\abs{b_m},  & \text{ if } m=n+1, \\
\frac{\abs{a_n b_{m-1}-a_{n-1}  b_{m}}} {(a_n \bar{a}_n +b_m \bar{b}_m)},  &  \text{ if } \text{m=n}.
\end{cases}
$$
\end{itemize}
\end{corollary}

It is worth to remark that the spherical derivative map is a bounded map since the domain is a compact topological space. Therefore, there is a suitable implementation of the spherical derivative map without overflows.

\begin{definition} Let $C=(x_0, f(x_0), \cdots, f^{q-1}(x_0))$ be a $q$-cycle of a rational map $f \colon \hat \C \to \hat \C$. The spherical derivative of $f$ at $C$ is defined as the product $$f^\# (C):=f^\# (x_0) f^\# (f(x_0)) \cdots  f^\# (f^{q-1}(x_0)).$$
A cycle $C$ is said to be \emph{super-attracting, attracting, indifferent or repelling} if  $f^\# (C)=0$, $f^\# (C)<1$, $f^\# (C)=1$ or $f^\# (C)>1$, respectively.
\end{definition}

%
\bigskip

 \subsection{Julia set and basins in a discrete semi-flow} \label{julia}

Let $(X,g)$ be a discrete semi-flow (shortly denoted as $X$) with (pseudo) metric $d$. The \emph{end space} of $(X,g)$ is the quotient
$\Pi(X,g)={X}/{\sim},$ where for $x,x'\in X$, $x\sim x' $ if, and only if,
$$\lim_{n\to+\infty} d(g^n(x),g^n(x'))=0.$$
An element $a=[x]\in\Pi(X,g)$ is called \emph{end point} of $X$. There is a canonical quotient map
\emph{$\omega_g\colon X\to\Pi(X,g),$} $\omega_g(x)=[x]$. The map $\omega_g$ gives rise to a decomposition of the discrete semi-flow into disjoint subsets
$$X=\bigsqcup_{a\in\Pi(X,g)}\omega_g^{-1}(a).$$
The subset $\omega_g^{-1}(a)$ with $a\in\Pi(X,g),$ is called the \emph{basin of the end point} $a$. In particular if $x_0 \in X$ is a fixed point, we have that the basin of $[x_0]$ is given by $$\omega_g^{-1}([x_0] )=\{x\in X\hspace{2pt} |\hspace{2pt} \lim_{n\to +\infty}d(g^n(x),x_0)=0\}.$$

\begin{definition} Let $(X,g)$ be discrete semi-flow and $T\subset X$ a subset. The \emph{basin determined by $T$} is
defined as
$$B_g(T):=\omega_g^{-1}(\omega_g(T)).$$
\end {definition}

We recall the notions of Fatou and Julia sets:

\begin{definition} Given  a discrete semi-flow $(X,g)$, a point $x_0 \in X$ is said to be a \emph{Fatou} point if there exists a neighbourhood $U$ at $x_0$ such that for every $\varepsilon >0$ there is $n_\varepsilon \in \N$ with
$$d(g^n(x), g^n(y)) <\varepsilon$$
\noindent for all $x, y \in U$ and $n\geq n_{\varepsilon }$. Denote as $F(g)$ the set of all Fatou points in $X$ and $J(g)=X\setminus F(g)$ its complement, usually called \emph{Julia} set of $g$.
\end {definition}

Given a point $x_0\in X$ the subset $\{g^n(x_0)\hspace{2pt}|\hspace{2pt} n \in \mathbb{N}\}$ is called the \emph{orbit} of $x_0$. Sometimes we also call orbit to the map $ \mathbb{N} \to X,$ $n\mapsto g^n (x_0)$. The \emph{backward orbit} of $x_0$ is the subset $\cup_{n \in \N} g^{-n}(x_0)$. For a given subset $S\subset X$ we may also consider its orbit $\cup_{n \in \N} g^{n}(S)$ and its backward orbit
$\cup_{n \in \N} g^{-n}(S)$.

For more properties about the Julia set associated to a rational map on the Riemann sphere we refer the reader to \cite{B1991},\cite{Milnor}.


\subsection{Ergodic measures on the Julia set of a rational map}

A probability measure $\mu$ supported on the Julia set $J(f)$ of a rational map $f:\hat\C\rightarrow \hat\C$ is called \emph{$\delta$-conformal} ($0<\delta\leq 2$) if
$$\mu(f(A))=\int_A \abs{f^{\#}}^\delta d\mu$$
\noindent holds for every Borel set $A \subset J(f)$ in which $f|_A$ is injective.

Sullivan proved in \cite[Theorem 3]{Sull83} the existence of a $\delta$-conformal measure:


\bigskip

{\it {\it Theorem 3}. There is a positive, finite measure $\mu$  on the Julia set
$J(f)$ of a rational map $z \mapsto f(z)$ satisfying for some real number $\delta=\delta(f)$
$$\mu(f(A))=\int_A \abs{f^{\#}}^\delta d\mu$$
for any Borel set $A\subset J(f)$ where $f$ is inyective. Moreover, $0< \delta(f)\leq 2 $ and if any $\delta$ satisfies the above equation for some measure $\mu$,  then $\delta \geq \delta (f)$.}

\bigskip

Now we recall the notion of an expanding rational map:

\begin{definition} Let
$f\colon \hat{\mathbb{C}} \to \hat{\mathbb{C}}$ be  a rational map.
We say that $f$ is \emph{expanding} (on $J(f)$) whenever there exists real numbers
$c>0,\lambda>1$ such that
$$(f^n)^{\#}(u)>c \lambda^n$$ \noindent for all $u \in J(f)$ and $1\leq n \in \N$.
\end{definition}

In the case of expanding rational maps, there is only one conformal ergodic measure
$ \mu$  which related to the Hausdorff dimension and Hausdorff measures as follows (see \cite[Theorem 4]{Sull83}):

\medskip
{\it {\it Theorem 4}. In the expanding case, there is one, and only one, conformal measure $\mu$ on the Julia set $J(f)$. The exponent $\delta=\delta(f)$ of $\mu$ is the Hausdorff dimension and $\mu$ is, up to a product with a constant, the Hausdorff $\delta$-measure $H_\delta$ on $J(f)$. Moreover, $0 < \delta(f) < 2$.}

\medskip
There are important results on the existence of ergodic measures on the Julia set of a rational map. For more results about these problems and related topics we may refer the reader, among others, to \cite{Sull83} and \cite{Graczyk2008NonuniformHI}.

\section{Lyapunov functions and basins}

Let $(X,f)$ be a discrete semi-flow and $\phi \colon X \to \R_+$ a continuous map. We can consider the functions given in subsection \ref{LyaTimeSpaceAverage}: $l_f(\phi)$,  $L_f(\phi )$, $l_f^q(\phi)$, whose domains are $\text{dom}(l_f(\phi ))=\{x \in X\hspace{2pt} |\hspace{2pt} l_f^-(\phi)(x)= l_f^+(\phi)(x)\}$ and similarly for $L_f(\phi)$ and $l_f^q(\phi)$. Also remember from subsection \ref{julia} the canonical quotient map
$\omega_f\colon X\to\Pi(X,f)$ from $X$ to the corresponding space of end points.

\begin{proposition}\label{fixedpoint} Let $(X,f)$ be a discrete semi-flow, $\phi \colon X \to \R_+$ a continuous map and $x_0$ a fixed point of $f$. The following statements hold true:

\begin{itemize}
\item[(i)] $\omega_f^{-1}(\omega_f (x_0) )\subset \text{dom}(l_f(\phi ))$ and $l_f(\phi )(x)=\phi(x_0)$, for all $x \in  \omega_f^{-1}(\omega_f (x_0) ).$

\item[(ii)] If $x_0 \in \phi^{-1}(0)$, then $$\omega_f^{-1}(\omega_f (x_0) )\subset \text{dom}(L_f(\phi ))
\cap  \text{dom}(l_f(\phi ))$$ \noindent and $l_f(\phi )(x)=L_f(\phi )(x)=\phi(x_0)=0$ for all $x \in  \omega_f^{-1}(\omega_f (x_0)).$

\item [(iii)] If $x_0 \not \in \phi^{-1}(0)$, then $$\omega_f^{-1}(\omega_f (x_0) ) \setminus  {\rm{linv}} (\phi^{-1}(0)) \subset \text{dom}(L_f(\phi ))
  \cap  \text{dom}(l_f(\phi ))$$ \noindent and
$l_f(\phi )(x)=L_f(\phi )(x)=\phi(x_0)$, for all $x \in  \omega_f^{-1}(\omega_f (x_0) ) \setminus  {\rm{linv}} (\phi^{-1}(0)).$
\end{itemize}
\end{proposition}

\begin{proof} If $x \in  \omega_f^{-1}(\omega_f (x_0) )$, then $\lim_n d(f^n(x), f^n(x_0))=0$ and therefore we have $\lim_n(f^n(x))= x_0$. Since $\phi$ is a continuous map, we have $\lim_n\phi(f^n(x))= \phi(x_0)$. Now (i), (ii) and (iii) follow from basic properties of limits and Corollary \ref{prod1}.
\end{proof}

\bigskip

We recall that for the discrete semi-flows $(X,f)$ and $(X,f^q)$ one has the corresponding end point spaces $\Pi(X,f)$ and $\Pi(X,f^q)$ as well as their respective canonical maps $\omega_f \colon X \to \Pi(X,f)$ and $\omega_{f^q} \colon X \to \Pi(X,f^q)$. It is straightforward to check that if
$(x_0, f(x_0), \cdots,  f^{q-1}(x_0))$ is a $q$-cycle of $f$  ($f^{q-1}(x_0)\not =x_0$ and $f^{q}(x_0)=x_0$), then $\omega_f^{-1}(\omega_f(x_0))=
\omega_{f^q}^{-1}(\omega_{f^q}(x_0))$. For any given $q$-cycle $C=(x_0, f(x_0), \cdots f^{q-1}(x_0))$, we consider the basin of $C$:
$$
\begin{array}{cll}
B_f(C) & := & \omega_{f}^{-1}(\omega_{f} (C)) \\
& = &\omega_{f}^{-1}(\omega_{f} (x_0)) \cup \omega_{f}^{-1}(\omega_{f}  (f(x_0))) \cup\cdots  \cup \omega_{f}^{-1}(\omega_{f}(f^{q-1}(x_0))).
\end{array}$$

\begin{proposition} Let $(X,f)$ be a discrete semi-flow, $\phi \colon X \to \R_+$  a continuous map and $C=(x_0, f(x_0), \cdots f^{q-1}(x_0))$ a $q$-cycle of $f$. The following statements hold true:

\begin{itemize}
\item[(i)] $B_f(C) \subset  \text{dom}(l_f^q(\phi ))$ and for every  $x\in B_f(C)$,
$$l_f^q(\phi) (x)= \big (\phi(x_0) \phi (f(x_0)) \cdots \phi (f^{q-1}(x_0)) \big )^\frac{1}{q}.$$

\item[(ii)] If $C \cap  \phi^{-1}(0)\not =\emptyset$, then $B_f(C) \subset \text{dom}(L_f(\phi ))
  \cap  \text{dom}(l_f^q(\phi ))$. Moreover, $l_f^q(\phi )(x)=L_f(\phi )(x)=0$, for all $x \in  B_f(C).$

\item [(iii)] If $C \cap  \phi^{-1}(0)=\emptyset$, then $B_f(C)  \setminus  {\rm{linv}} (\phi^{-1}(0)) \subset \text{dom}(L_f(\phi ))
\cap  \text{dom}(l_f^q(\phi ))$. Moreover,
$l_f^q(\phi )(x)=L_f(\phi )(x)=  \big (\phi(x_0) \phi (f(x_0)) \cdots \phi (f^{q-1}(x_0)) \big )^\frac{1}{q}$
for all $x \in  B_f(C)  \setminus  {\rm{linv}} (\phi^{-1}(0)).$

\end{itemize}

%
\end{proposition}
\begin{proof} This is just a consequence of Corollary \ref{arith5} .
\end{proof}

Now for a discrete semi-flow $(X,f)$ and a  continuous  map $\varphi \colon X \to \{-\infty\}\cup  \R$, consider the corresponding functions $a_f(\varphi)$, $A_f (\varphi)$, $a_f^q (\varphi)$ taking values in $\{-\infty\}\cup  \R \cup \{+\infty\}$ (see subsection \ref{LyaTimeSpaceAverage}).

\begin{proposition} Let $(X,f)$ a discrete semi-flow, $\varphi \colon X \to \{-\infty\}\cup  \R$ a continuous map and $C=(x_0, f(x_0), \cdots f^{q-1}(x_0))$ a $q$-cycle of $f$. Then,

\begin{itemize}
\item[(i)] $B_f(C)\subset  \text{dom}(a_f^q(\varphi ))$ and
$$a_f^q(\varphi) (x)= \frac{\varphi (x_0)+ \varphi (f(x_0))+ \cdots +\varphi(f^{q-1}(x_0))}{q}$$
\noindent for all $x\in B_f(C)$.

\item[(ii)]  If $C \cap  \varphi^{-1}(-\infty)\not = \emptyset$, then $ B_f(C) \subset \text{dom}(A_f(\varphi ))
\cap  \text{dom}(a_f^q(\varphi ))$ and  $$a_f^q(\varphi )(x)=A_f(\varphi )(x)=-\infty$$ \noindent for all $x \in  B_f(C)$.

\item [(iii)] If $C \cap  \varphi^{-1}(-\infty)=\emptyset$, then $B_f(C)\setminus  {\rm{linv}}(\varphi^{-1}(-\infty))\subset \text{dom}(A_f(\phi ))
\cap  \text{dom}(a_f^q(\phi ))$ and $$a_f^q(\phi )(x)=A_f(\phi )(x)=  \frac{\varphi (x_0)+ \varphi (f(x_0))+ \cdots +\varphi(f^{q-1}(x_0))}{q}$$
\noindent for all $x \in  B_f(C)  \setminus  {\rm{linv}} (\varphi^{-1}(\infty))$.

\end{itemize}

\end{proposition}

Now suppose that we consider the Borel $\sigma$-algebra $\Sigma$ on $X$ for a given discrete semi-flow $(X,f)$. If $B_f(C) =\omega_f^{-1}(\omega_f (C) ) \in \Sigma$ for any $q$-cycle $C$, then $B_f(C)$ is a discrete sub-semi-flow with a continuous self-map $f|_{B_f(C)} \colon B_f(C) \to B_f(C)$. Moreover,  there is an induced Borel $\sigma$-algebra on $B_f(C)$.
If $C=(x_0, f(x_0), \cdots, f^{q-1}(x_0))$ stands for a $q$-cycle of $f$, then we can take the discrete Borel measure $\mu _C=\sum _{k=0}^{q-1}\delta _{f^k(x_0)}$ on $B_f(C)$ and the Space Average
$$S(\varphi|_{B_f(C)}):=\frac {1}{q}\int_{B_f(C)} \varphi \,d\mu_C=  \frac {\varphi(x_0)+\varphi(f(x_0))+ \cdots+ \varphi(f^{q-1}(x_0))}{q}$$
\noindent for a continuous map $\varphi \colon X \to  \{-\infty\}\cup  \R$.

\begin{proposition} Under the conditions and notation described above the following statements hold:
\begin{itemize}
\item[(i)] $B_f(C) \setminus (\bigcup_{k>0} f^{-k}(\varphi^{-1}(-\infty)) ) \subset \text{dom}(A_f(\varphi))$ and
$$A_f(\varphi)(x)=\lim_{n \rightarrow + \infty}\; \frac{1}{n} \sum_{k=0}^{n-1} \varphi\left(f^k (x)\right) =  \frac {\varphi(x_0)+\varphi(f(x_0))+ \cdots+ \varphi(f^{q-1}(x_0))}{q}$$ \noindent for  all $x \in B_f(C)  \setminus (\bigcup_{k>0} f^{-k}(\varphi^{-1}(-\infty)) )$.

\item[(ii)] $A(\varphi)(x)=S(\varphi|_{B_f(C)})$ for  all $x \in B_f(C) \setminus (\bigcup_{k>0} f^{-k}(\varphi^{-1}(-\infty)) )$. In other words, $$\lim_{n\rightarrow +\infty}\; \frac{1}{n} \sum_{k=0}^{n-1} \varphi\left(f^k (x)\right) =\frac {1}{q}\int_{B_f(C)} \varphi \,d\mu_C.$$
\end{itemize}

\end{proposition}

As an immediate consequence we have:

\begin{corollary} 
Under the conditions and notation described above, $$(B_f(C),\Sigma, f, \mu_C)$$\noindent is a $B$-space.
\end{corollary}

\begin{remark}
If the continuous map $\phi \colon X \to [0,+\infty)$ satisfies the finite-nullity condition, then
$$\lim_{n\rightarrow +\infty}\; \frac{1}{n} \sum_{k=0}^{n-1} \log \phi \left(f^k (x)\right) = \frac { \log \phi (x_0)+\log \phi (f(x_0))+ \cdots+ \log \phi (f^{q-1}(x_0))}{q}$$
$$= \frac {1}{q}\int_{B_f(C)} \log \phi \,d\mu_C $$
\noindent for all $x\in B_f(C)\setminus (\bigcup_{k>0} f^{-k}(\phi^{-1}(0)) )$.

It is interesting to note that is $X$ is a metric space, under this condition, $\bigcup_{k>0} f^{-k}(\phi^{-1}(0))$ is an isolated countable subspace of $X$. Observe that, if $X=S^2$, then every isolated countable subspace has Lebesgue measure zero.
\end{remark}

\begin{remark}
If $f\colon S^{2} \to  S^{2}$ is a rational map such that $f|_{J(f)}$ is ergodic with respect to a measure $\mu$ defined on
$J(f)$, then $S^{2}$ can be expressed as a finite sum (i.e. disjoint union) of $B$-spaces $$S^{2}=(J(f),\mu) +(B_f(C_1), \mu_{C_1}) +\cdots +(B_f(C_m), \mu_{C_m})$$ \noindent
where each $C_i$ is an attracting cycle and $\mu_{C_i}$ is a discrete measure (see Remark \ref{quasi}).
\end{remark}

\section{Computing basins of attraction  of $q$-cycles, Lyapunov constants and Lyapunov exponents of a rational map}\label{algorithms}

Taking the complex projective line $X=\p^1(\C)$ as a model, and an r-homogeneous irreducible pair $(F,G) \in \C_h[z,t]\underline{\times}_{\rm{ri}} \C_h[z,t]$ (see Definitions  \ref{bivariate}, \ref{homogeneouspair} and \ref{irreducible}) we know that there are induced maps (see \eqref{erre}, \eqref{ese} and \eqref{pe}):
$$\begin{array}{l}
R\colon \C^2 \setminus \{(0,0)\} \to \C^2 \setminus \{(0,0)\},\\
R^S\colon S^{3}\to S^{3}, \\
R^P\colon \p^{1}(\C)\to \p^{1}(\C).
\end{array}$$

In turn, associated to $f=R^P\colon \p^{1}(\C)\to \p^{1}(\C)$ and a continuous map $\phi \colon \p^{1}(\C) \to \R_+$ we can consider the functions $l_{R^P}(\phi)$,  $L_{R^P}(\phi )$, $l_{R^P}^q(\phi)$ given in subsection \ref{LyaTimeSpaceAverage}. We note that the quotient maps
$$\xymatrix{\C^2 \setminus \{(0,0)\} \ar[r]^(.7)p & S^{3} \ar[r]^(.4)q & \p^{1}(\C)}$$
\noindent induce the composite $\phi^S=\phi \circ q \colon S^{3} \to \p^{1}(\C) \to \R_+$.
Now, if $r, s \in \mathbb{N}$ with $r\leq s,$ then we can also take the map $L_{R^P}^{[r,s]}(\phi) \colon  \p^{1}(\C) \to \mathbb{R}_+$ given as
$$L_{R^P}^{[r,s]}(\phi) (x)=\left(\prod_{k=r}^s \phi ((R^P)^k(x)) \right )^{\frac{1}{s+1-r}}, \quad x \in  \p^{1}(\C)$$
\noindent as well as the induced map $L_{R^S}^{[r,s]}(\phi^S) \colon  S^3 \to \mathbb{R}_+$, given as
$$L_{R^S}^{[r,s]}(\phi^S) (y)=\left (\prod_{k=r}^s \phi^S ((R^S)^k(q(y))) \right )^{\frac{1}{s+1-r}}, \quad y \in S^3. $$

Let $\mathcal{C}=\{C_1, \cdots C_m\}$ denote the finite set of attracting cycles of $f$. Given a cycle $C_i= (x_i,f(x_i)\cdots,f^{q_i-1}(x_i))$ we obtain what we call the \emph{Lyapunov constants} and \emph{Lyapunov exponents}, respectively: $$c_i=(\phi (x_i)\phi (f(x_i))\cdots \phi (f^{q_i-1}(x_i)))^{\frac{1}{q_i}},$$
$$\log (c_i)=\frac{\log \phi (x_i)+ \log \phi (f(x_i))+\cdots+ \log  \phi (f^{q_i-1}(x_i))}{q_i} .$$ We may clearly assume that $c_1 \leq c_2 \leq \cdots \leq  c_m$.  We want to point out that
$$L_f(\phi)|_{B_f(C_i)}=c_i$$
\noindent almost everywhere. Now, assume $1\leq m_1< \cdots <m_r=m$ are such that $$c_1=\cdots=c_{m_1}<c_{m_1+1}=c_{m_1+2}=\cdots=c_{m_2}< \cdots< c_{m_{r-1}+1}=\cdots= c_{m_r}=c_m.$$
\noindent Then, up to backward orbit of critical points, one has:
$$L_f(\phi )({\cup_{i=1}^{m_1}B_f(C_i)})=c_{m_1},$$
$$L_f(\phi )({\cup_{i=m_1+1}^{m_2}B_f(C_i)})=c_{m_2},$$
$$\vdots$$
$$L_f(\phi )({\cup_{i=m_{r-1}+1}^{m_r}B_f(C_i)})=c_{m_r}.$$

If $q(y)\in \cup_{i=m_j+1}^{m_{j+1}}B_f(C_i)$, then,  for every $C \in \N$ there exists $k_{C}\in \N$ such that
\begin{equation}
\label{eq2}
\abs{L_{f}^{[0,k-1]}(\phi^S)(y)-L_{f}^{[0,k]}(\phi^S)(y)}\leq 10^{-C},\mbox{ for every }k\geq k_{C}.
\end{equation}

And using the quotient map $q\colon S^3 \to \p^1(C)\cong S^2$, we can take the inverse image $q^{-1}(\cup_{k>0} f^{-k}(\phi^{-1}(0))).$ We observe that
$$L_{f^S}(\phi^S) (y)= \lim_{s\to +\infty}  L_{f^S}^{[r,s]}(\phi^S)(y) \hspace{5pt}(r\in  \N)$$
\noindent holds whenever $y \in (S^3 \setminus (q^{-1}(\cup_{k>0} f^{-k}(\phi^{-1}(0))))) \cap \text{dom} L_{f^S}(\phi^S)$.
Therefore given $r$ and $C \in \N$, one can find $k_{(C,r)}$ such that

\begin{equation}
\label{eq2r}
\abs{L_{f^S}^{[r,k-1]}(\phi^S)(y)-L_{f^S}^{[r,k]}(\phi^S)(y)}\leq 10^{-C},\mbox{ for every }k\geq k_{(C,r)}.
\end{equation}
Also, if $y \in q^{-1} (B_f (C_i)\setminus (\cup_{k>0} f^{-k}(\phi^{-1}(0))))$ and $C_i=(x_i, \cdots, f^{q_i-1} (x_i)) $ is a $q_i$-cycle, then
$$L_{f^S}(\phi^S) (y)=  l_{f^S}^{q_i}(\phi^S)(y) = \lim_{s\to +\infty} L_{f^S}^{[s-1,s-q_i]}(\phi^S)(y). $$

Consequently, given $C \in \N$, one can find $k_{C}$ such that
\begin{equation}
\label{eq2r}
\abs{L_{f^S}^{[k-1,k-q_i]}(\phi^S)(y)-L_{f^S}^{[k,k+1-q_i]}(\phi^S)(y)}\leq 10^{-C}, \mbox{ for every }k\geq k_C.
\end{equation}

Then, an algorithm that computes the basins can be developed as we will describe in the next subsection.

\subsection{Computing Lyapunov constants and $q$-cycles}

This first part of our method can be done by considering a representative finite ``grid'' $D \subset \p^{1}(\C)$. For each $x\in D$ ($x=[z:t],\, \abs{z}+\abs{t}=1)$ we can take a point $y \in S^3$ ($y=(z,t)$) such that $q(y)=x$. Consider this finite set $D^S\subset S^3$.

We can cover $\p^{1}(\C)$ by using two rectangular parametrizations. We contemplate one representative grid of the form $[z_{ij}:1]$ and another of the form $[1:t_{ij}]$. Considering normalized coordinates we obtain points $$(\frac{z_{ij}}{1+\abs{z_{ij}}}, \frac{1}{1+\abs{z_{ij}}}),\hspace{6pt} (\frac{1}{1+\abs{t_{ij}}}, \frac{t_{ij}}{1+\abs{t_{ij}}})$$ \noindent in $S^3$. Since $\p^{1}(\C)\cong S^2$, taking $S^2$ as a model one can regard cubic  or simplicial structures and subdivision procedures to construct a starting ``grid'' of points in $S^2$. This second procedure of selecting a starting grid seems to be a more suitable and ``representative'' grid of points in $S^2$.

Now, given a point $y \in D^S$ there are two possible cases: If $L_{f^S}(\phi^S)(y)$ is defined, then we take the constant $c=L_{f^S}(\phi^S)(y)$; otherwise, we try with another point $y'\in D^S$.

We know that if $L_{f^S}(\phi^S)(y)$ is defined, then one has that, up to some tolerance $10^{-C}$, there exists $k>>0$ (i.e., much greater than $0$) such that
$$\abs{L_{f^S}^{[0,k-1]}(\phi^S)(y)-L_{f^S}^{[0,k]}(\phi^S)(y)} \leq 10^{-C}.$$
In the case that the finite-nullity condition (see Definition  \ref{fphicondition}) is verified, for any $r\in \N$ we can check whether there is $k>r$ such that
$$\abs{L_{f^S}^{[r,k-1]}(\phi^S)(y)-L_{f^S}^{[r,k]}(\phi^S)(y)} \leq 10^{-C}.$$

The second part of this method is devoted to developing algorithms that divide the complex projective line into $1+ r$ disjoint regions induced by the constants $c_{m_1}<\cdots<c_{m_r}$. Furthermore, these algorithms are  refined to give a partition the sphere into $m+1$ regions associated with the corresponding attracting  cycles  $\mathcal{C}=\{C_1, \cdots C_m\}$ and the  complementary of the union of the basins of attraction.

\subsubsection{Lyapunov constants associated to fixed points}

Given a tolerance $10^{-C}$, a maximum number of iterations $M$ and an initial subset of starting points $D^S$ we follow the next procedure:

For a given point $y \in D^S$, one can check whether there is $k_{(C,M)}< M$ satisfying the following condition:
$$\abs{\phi^S((f^S)^{k_{(C,M)}-1}(y))- \phi^S((f^S)^{k_{(C,M)}}(y))}\leq 10^{-C} \text{ and }$$
\vspace{-25pt}\begin{equation}
\label{Cauchyphiand}
\abs{\phi^S((f^S)^{k-1}(y))- \phi^S((f^S)^{k}(y))}> 10^{-C}, \quad 0<k< k_{(C,M)}.
\end{equation}
This condition will be referred to as the \emph{$M$-truncated Cauchy inequality} for $\phi^S$.

Depending on whether condition \eqref{Cauchyphiand} is satisfied or not, we have two cases:
\begin{enumerate}
\item[a)] If there is $k=k_{(C,M)}<M $, then the output is the pair $((f^S)^k(y),k)$ and we take $(f^S)^k(y)$ as a member of a new list of fixed points (1-cycles) $\text{Cycles}_1$ and we also take the constant $c^1=\phi^S((f^S)^k(y)).$

\item[b)] Otherwise, we take $k_{(C,M)}=M $ and we construct a new list $D_1^S$ and $(f^S)^M(y)$ is a new member of the list $D_1^S$.
\end{enumerate}
When $y$ has run over all the points in $D^S$ we will have constructed the lists $D_1^S$, $\text{Cycles}_1$ and $c_1^1, \cdots, c_{i_1}^1$.

\subsubsection{Lyapunov constants associated to $q+1$ cycles} \label{EvaluationCauchy}

Suppose that we have computed the lists below:

\medskip
$D_1^S$, $\text{Cycles}_1$ and $c_1^1, \cdots, c_{i_1}^1$ (1-cycles)

$D_2^S$, $\text{Cycles}_2$ and $c_1^2, \cdots, c_{i_2}^2$ (2-cycles)

$ \text{                   }  \vdots $

$D_q^S$, $\text{Cycles}_q$ and $c_1^q, \cdots, c_{i_q}^q$ ($q$-cycles).

\medskip
For a given point $y \in D_q^S$, we can check whether there is $k_{(C,M)}<M$ satisfying the following condition:
%
\begin{equation*}
\label{}
\abs{L_{f^S}^{[k_{(C,M})-q-1,k_{(C,M})-1]}(\phi^S)(y)-L_{f^S}^{[k_{(C,M)}-q,k_{(C,M)}]}(\phi^S)(y)} \leq 10^{-C} \text{ and }
\end{equation*}
\vspace{-15pt}\begin{equation}
\label{Cauchyphiand1}
\abs{L_{f^S}^{[k-q-1,k-1]}(\phi^S)(y)-L_{f^S}^{[k-q,k]}(\phi^S)(y)} > 10^{-C}, \quad q\leq k< k_{(C,M)}.
\end{equation}

There are two cases:
\begin{enumerate}
\item[a)] If the condition above is satisfied for $k_{(C,M)}<M $, then the output is the pair $(((f^S)^{k_{(C,M)}-q}(y), \cdots,  (f^S)^{k_{(C,M)}}(y)),k_{(C,M)})$ and we take $$((f^S)^{k_{(C,M)}-q}(y), \cdots,  (f^S)^{k_{(C,M)}}(y))$$ \noindent as a member of a new list $\text{Cycles}_{q+1}$. We also have the constant $c^{q+1}=L_{f^S}^{[k_{(C,M)}-q,k_{(C,M)}]}(\phi^S)(y).$

\item[b)] Otherwise, take $k_{(C,M)}=M $ and we construct a new list $D_{q+1}^S$, $(f^S)^M(y)$ being a new member of the list $D_{q+1}^S$.
\end{enumerate}

When $y$ has run over all the points in $D_q^S$ we will have constructed the lists $D_{q+1}^S$, $\text{Cycles}_{q+1}$ and $c_1^{q+1}, \cdots, c_{i_{q+1}}^{q+1}$.

\bigskip

It is plain to check that for a rectangle $(a,b) \times (c,d)$ containing the unit disk, the images of the parametrizations
$$\begin{array}{l}
\alpha \colon (a,b) \times (c,d)\to  \p^1(\C),\hspace{6pt} (u,v)\mapsto [u+ i v : 1] \\
 \beta \colon (a,b) \times (c,d)\to  \p^1(\C), \hspace{6pt} (u,v)\mapsto [1:u+ i v ]
\end{array}$$ \noindent cover $\p^1(\C)$, and we can choose as initial starting points of the form
$$(\frac{u+ i v}{1+(u^2+v^2)^{1/2}}, \frac{1}{1+(u^2+v^2)^{1/2}})\hspace{6pt}\mbox{or}\hspace{6pt}(\frac{1}{1+(u^2+v^2)^{1/2}}, \frac{u+ i v}{1+(u^2+v^2)^{1/2}}).$$

We want to mention that we may need to refine the lists $\text{Cycles}_{s}$ to new ones $\text{Cycles}'_{s}$ as we may have different $s$-cycles being equal up to tolerance $10^{-C}$. This produces the list of cycles
$\text{Cycles}'_{1}\cup \text{Cycles}'_{2} \cup \cdots \cup \text{Cycles}'_{m}$.

Note that using the spherical derivative we know whether an $s$-cycle is  attracting. Then we can have a second reduction of these lists and take attracting $s$-cycles only.

\begin{remark} Other techniques based on spherical multipliers (related to geometrical means) 
 that,  instead of computing basins of attracting cycles, give a graphic representation of neighbourhoods of the Julia set have been  used by the authors in
\cite{EXTREMIANAALDANA2018251} and \cite{ExtremianaAldana2018SelfoverlaysAS} . Our new methods, based on computing Lyapunov functions, have been implemented in Julia language, see  \cite{Alvarez-Aparicio_Lyapunov_Cycle_Detector_2022}. We remark that  some basic  functions and procedures used  to write the new code  had been  previously developed in \cite{hernandez19},  \cite{EGHR19}.
\end{remark}

\section{Applications to Newton's method}

Newton's method, besides being one of the most used numerical methods to approximate solutions to non-linear equations, it is also a method that has been profusely studied, both from a numerical and a dynamic perspective. In this regard, many useful results can be found in the literature (\cite{Sull83},\cite{Var02}, \cite{Drexler}, and many more). As consequence,  we think that it is interesting to give a new  version which consider  the Newton's method  as a rational Hopf-endomorphism.  We remark that in the case of a polynomial with simple roots, the homogeneous pair induced by Newton's method is always $r$-irreducible and we can apply the results developed in this paper.  For  other results  about the homogeneous Newton's method or for   the case of multiple roots,  see \cite{GGHR15}, \cite{HernndezParicio2015BivariateNM}, \cite{EGHR2021}.


\subsection{The homogeneous Newton's method}\label{rhN}

Given a polynomial $p\in \C[x]$, $p(x)=a_0+a_1x+\cdots +a_n x^n$, $a_n\not =0$,
its standard derivative polynomial $p'$ is $p'(x)=a_1+ 2 a_2 x+\cdots +n a_n x^{n-1}$.
We can consider what we call the \emph{Homogeneous Newton's method}: $$\text{HN}: \C[x] \to \C_h[z,t]\underline{\times} \C_h[z,t]
,\hspace{10pt}\text{HN}(p):=({\rm H}(x p'- p,d) , {\rm H}(p', d))$$ \noindent where $d=\text{degree}(p)$ and
$\text{ H}: \C[x] \underline{\times} \N \to \C_h[z,t]$ is the homogenization map defined in equation \eqref{homogeneousversion}.
In this case, we obtain the following homogeneous pair of degree $n$:
\begin{equation}
\label{ }
\text{HN}(a_0+a_1x+\cdots +a_n x^n)=(F(z,t),G(z,t)).
\end{equation}
\vspace{-25pt}
\begin{eqnarray}
F(z,t)&=&-a_0 t^n+(a_1-a_1) z t^{n-1}+\cdots +(na_n- a_n) z^n, \\
G(z,t)&=&a_1 t^n+ 2 a_2 z t^{n-1}+ \cdots + n a_n z^{n-1}t.
\end{eqnarray}

\medskip
Notice that if $p=a_0+a_1x+\cdots+a_n x^n$ ($a_n\not=0$) is a polynomial where all its roots are simple and $\text{degree}(p)\geq 1$, then
the homogeneous pair induced by Newton's method 
 $\text{HN}(p)=({\rm H}(x p'- p,d),{\rm H}(p', d))$ is $r$-irreducible. This is a case we are specially interested in.
 
In order to work with a normalized homogeneous pair we shall consider the pair $R=(\bar F(z,t),\bar G(z,t)$:
\begin{eqnarray}
\label{FGnormalized}
\bar F(z,t)=\frac{-a_0 t^n+(a_1-a_1) z t^{n-1}+\cdots +(na_n- a_n) z^n}{\abs{-a_0}+ \abs{a_1-a_1} + \cdots +\abs{na_n- a_n} + \abs{a_1}+ \abs{2 a_2}+ \cdots +\abs{n a_n }},
\end{eqnarray}
\begin{eqnarray}
\label{FGnormalized1}
\bar G(z,t)=\frac{a_1 t^n+ 2 a_2 z t^{n-1}+ \cdots + n a_n z^{n-1} t }{\abs{-a_0}+ \abs{a_1-a_1} + \cdots +\abs{na_n- a_n} + \abs{a_1}+ \abs{2 a_2}+ \cdots +\abs{n a_n }}.
\end{eqnarray}

Now we have the following self-map of the Hopf fibration
$$\xymatrix{
\ar[d]_q S^{3} \ar[r]^{R^S}  & S^{3} \ar[d]^q \\
\p^{1}(\C) \ar[r]^{R^P} &\p^{1}(\C) }
$$
\noindent where, taking into account equations \eqref{FGnormalized} and \eqref{FGnormalized1}, $R^S$ and $R^P$ are given as
\begin{eqnarray}
\label{selfmaphopf}
R^S((z,t))=(\frac{\bar F(z,t)}{\abs{\bar F(z,t) }+\abs{\bar
G(z,t) }},   \frac{\bar G(z,t)}{\abs{\bar F(z,t) }+\abs{\bar
G(z,t) }}),
\end{eqnarray}
\begin{eqnarray}
\label{selfmaphopf1}
 R^P([z:t]))=[\frac{\bar F(z,t)}{\abs{\bar F(z,t) }+\abs{\bar
G(z,t) }}:  \frac{\bar G(z,t)}{\abs{\bar F(z,t) }+\abs{\bar
G(z,t) }}].
\end{eqnarray}

\subsection{The spherical derivative for the homogeneous Newton's method}\label{sphericalderivativesdiagram}

In this case, for the standard Newton's method we have $f=N_p:\hat{\C } \rightarrow \hat{\C }$, where
$$N_p(x)=\frac{x p'(x)- p(x)}{p'(x)}.$$
Obviously, the numerator and denominator of this rational map are $A(x)= x p'(x)- p(x)$ and $B(x)=p'(x)$ respectively. Moreover, its standard derivative is
$$N'_p(x)=\frac{A'(x) B(x)- A(x) B'(x)}{(B(x))^2}=\frac{p(x) p''(x)}{(p'(x))^2}.$$

Using the notation of subsection \ref{computingsphericalderivative} we have that $d=\text{degree}(p)=n$ and $e=2d-1$ with $e^*=e-1.$ Furthermore,
$$\begin{array}{l}
U(z, t)=U_p(z,t) U_{p''}(z,t) (t \bar t+ z \bar z), \\
V(z, t)=V_A(z,t) {\bar V}_A(\bar z,\bar t)+V_B(z,t) {\bar V}_B(\bar z,\bar t)
\end{array}$$
\noindent where
$$\begin{array}{l}
U_p(z,t):=t^{n} p(z/t) = a_0 t^n+a_1z t^{n-1}+\cdots +a_n z^n, \\
U_{p''}(z,t):=t^{n-2}p''(z/t)=2 a_2  t^{n-2}+6a_3z t^{n-3}+\cdots +n (n-1) a_n z^{n-2}, \\
V_A(z,t):=-a_0 t^n+(a_1-a_1) z t^{n-1}+ (2 a_2- a_2) z^2 t^{n-2}+\cdots +(na_n- a_n) z^n, \\
V_B(z,t):=t^{n} p'(z/t) = a_1 t^{n}+ 2 a_2 z t^{n-1}+ 3 a_3 z^2 t^{n-2}+\cdots +n a_n z^{n-1} t.
\end{array}$$
Accordingly, for Newton's method we have

\begin{equation}
\label{definitionUV}
\abs{\frac{U(z,t)}{V(z,t)}}= \abs{\frac{U_p(z,t) U_{p''}(z,t) (t \bar t+ z \bar z)}{V_A(z,t) {\bar V}_A(\bar z,\bar t)+V_B(z,t) {\bar V}_B(\bar z,\bar t)}}.
\end{equation}

Observe that the following diagram is commutative
$$
\xymatrix{ \ar[d]_q S^3  \ar[rd]^{\phi^S} & \\
\p^1(\C) \ar[r]_{\phi^P}& \R_+}
$$
\noindent where $\phi^S(z,t)= \abs{\frac{U(z, t)}{V(z,t)}}$  and $\phi^P([z:t])= \abs{\frac{U(z, t)}{V(z,t)}}$.
Moreover, the spherical derivative $f^\#:\hat{\C}\rightarrow \mathbb{R}_+$ is given by
$$
f^\#(x)=\begin{cases}
\phi^S(\frac{x} {1+\abs{x}},\frac{1} {1+\abs{x}})(=\phi^P([\frac{x} {1+\abs{x}}:\frac{1} {1+\abs{x}}]), & x \in \C, \\
\frac{n}{n-1}, &  x=\infty.
\end{cases}
$$

Note that if $p(x_0)=0$ or $p''(x_0)=0$, then $f^\#(x_0)=0$. Certainly, for Newton's method applied to a polynomial $p$ with degree $n>1$, the roots of $p(x)=0$ and the roots of $p''(x)=0$ are critical points and $\infty$ is a repelling point.

\begin{example}\label{example1} If $\alpha \in \C,$ then consider the cubic polynomial $$p_{\alpha}(x)=(x-\alpha)(x^2+1)=-\alpha+x-\alpha x^2+x^3.$$
The homogeneous pair ${\rm{HN}}(p_{\alpha})$ associated to Newton's method is given by
\begin{eqnarray}
F_{\alpha}(z,t) & = &\alpha t^3- \alpha z^2 t+2 z^3, \\
G_{\alpha}(z,t) & = &   t^3- 2 \alpha z t^{2}+ 3 z^2 t.
\end{eqnarray}
$${\rm{HN}}(p_{\alpha})(z,t)=(F_{\alpha }(z,t) , G_{\alpha }(z,t))=(\alpha t^3- \alpha z^2 t+2 z^3,t^3- 2 \alpha z t^{2}+ 3 z^2 t ),\hspace{2pt}(z,t)\in S^3$$
\noindent and its normalization ${\rm{HN}}^S \colon S^3 \to S^3$ is identified with
\begin{equation}
\label{HNS }
{\rm{HN}}^S(p_{\alpha}) ((z,t))=(\frac{\alpha t^3- \alpha z^2 t+2 z^3}{4\abs{\alpha}+6} , \frac{ t^3- 2 \alpha z t^{2}+ 3 z^2 t  }{{4\abs{\alpha}+6}}), \quad (z,t)\in S^3.
\end{equation}
Finally, the homogeneous spherical derivative of the rational map is given by either $\phi_\alpha^S$ or $\phi_\alpha^P$, where
\begin{equation} \label{cubical}
\phi_\alpha^S(z,t)=\abs{\frac{U(z,t)}
{V(z,t)}}.
\end{equation}

$$\begin{array}{lcl}
U(z, t)&=&(- \alpha t^3+z t^{2}-\alpha z^2 t^{}+ z^3) (-2 \alpha  t^{}+ 6z) (t \bar t+ z \bar z), \\
V(z,t)&=&(\alpha t^3-\alpha z^2 t+2z^3)(\bar \alpha {\bar t}^3- \bar \alpha {\bar z}^2 {\bar t}+2{\bar z}^3)+ \\
      & & +(t^3- 2 \alpha z t^{2}+ 3 z^2)({\bar t}^{3}- 2 \bar \alpha {\bar z} {\bar t}^{2}+ 3  {\bar z}^2 {\bar t}).
\end{array}$$ 

Observe that $(\frac{i}{2} , \frac{1}{2})$, $(\frac{-i}{2} , \frac{1}{2})$, $(\frac{\alpha}{1 +\abs{\alpha}} , \frac{1}{1 +\abs{\alpha}})$  and $( \frac{\alpha}{3+\abs{\alpha}}, \frac{3}{3+\abs{\alpha}})$ are the  critical points. Note that the Newton's method associated to a cubic polynomial with simple roots $p(x)=(x-\alpha_1)(x-\alpha_2)(x-\alpha_3)$ has as free critical point the arithmetic mean $\frac{\alpha_1+\alpha_2+\alpha_3}{3}$ of the roots of $p(x)$. For $p_\alpha(x)$, one has  
$\frac{\alpha_1+\alpha_2+\alpha_3}{3}= \frac{\alpha+i-i}{3}= \frac{\alpha}{3}$ and  
$[\frac{\alpha}{3+\abs{\alpha}}: \frac{3}{3+\abs{\alpha}}]=[\frac{\alpha}{3}:1].$

\end{example}

\section{Basins of attraction  of $n$-cycles for the Newton's method associated to cubic polynomials with simple roots}\label{cubic_polynomials}

In this section we study the existence of attracting $n$-cycles for Newton's method applied to cubic polynomials.
We recall that Smale~\cite{Smale1985OnTE} introduced the concept of general convergence of an iterative method (later generalized by McMullen~\cite{McMullen87}), and proved that Newton's method is not generally convergent for cubic polynomials; i.e. Newton's method does not satisfy that for every initial point, and for almost every cubic polynomial, the method converges to a root of such polynomial.

In this section we focus in the case of cubic polynomials with single roots and we  aboard its study using rational self-maps of the Hopf fibration and Lyapunov exponents. 
As stated in \cite{ABP2004Scaling}, Newton's method satisfies the Scaling Theorem. As a consequence, our study can be essentially reduced to a one-parametric family of Newton's methods.

\subsection{The one-parametric family of Newton's method associated to cubic polynomials of the form $p_{\alpha}(x)=(x-\alpha) (x^2+1)$}

Given a complex cubic polynomial $p(x)$ with simple roots, one can always apply a M\"obius transformation to $p$ carrying two of its roots to the complex numbers $i$ and $-i$, the imaginary unit and its conjugate. The resulting polynomial is of the form $p_\alpha(x)=(x-\alpha)(x^2+1)$.
It is also important to bear in mind the following theorem established by Fatou~\cite{Fatou}

\medskip
 \emph{Theorem. Let $f:\hat{\C}\rightarrow \hat{\C}$ be a rational map. If $f$ has an attracting $n$-cycle, then the orbit of at least one critical point of $f$ converges to the $n$-cycle under the iteration of $f$.}

\medskip
We have seen in Example \ref{example1} that the normalized homogeneous coordinates of the critical points of Newton's method applied to $p_\alpha$ are given by the super-attractor fixed points
$[\frac{i}{2}:\frac{1}{2}]$, $[\frac{-i}{2}:\frac{1}{2}]$, $[\frac{\alpha}{1 +\abs{\alpha}}:\frac{1}{1 +\abs{\alpha}}]$ as well as the additional free critical point $[ \frac{\alpha}{3+\abs{\alpha}}:\frac{3}{3+\abs{\alpha}}]$. Our algorithms consider the free critical  point $y_\alpha = (\frac{\alpha}{3+\abs{\alpha}}, \frac{3}{3+\abs{\alpha}})\in S^3$ and the map $R_\alpha^S \colon S^3 \to S^3$ to construct the following orbit of  $y_\alpha\in S^3$:
$$(y_\alpha, R_\alpha^S(y_\alpha), (R_\alpha^S)^2(y_\alpha), \cdots, (R_\alpha^S)^k(y_\alpha), \cdots )$$
\noindent or the corresponding orbit of the point $q(y_\alpha)\in \p^1(\C)$:
\begin{equation}
\label{orbitcirtic}
(q(y_\alpha), R_\alpha^P(q(y_\alpha)), (R_\alpha^P)^2(q(y_\alpha)), \cdots, (R_\alpha^P)^k (q(y_\alpha)), \cdots ).
\end{equation}


Using the homogeneous spherical derivate $\phi_\alpha^S \colon S^3 \to \R_+ $ we consider the sequence
$$(\phi_\alpha^S(y_\alpha), \phi_\alpha^S(R_\alpha^S(y_\alpha)), \phi_\alpha^S((R_\alpha^S)^2(y_\alpha)), \cdots, \phi_\alpha^S((R_\alpha^S)^k(y_\alpha)), \cdots ).$$

Such algorithms also analyse the convergence of sequences of the form:

$$a_s^r=L_{R^S}^{[r,s]}(\phi^S) (R_\alpha^S(y_\alpha)), \quad \log a_s^r,\hspace{2pt}r>0.$$

The plane $\mathbb{C}$ of complex numbers  will be our parameter space in this setting ($\alpha \in \C$). In order to study the cyclic behaviour of Newton's method for cubic polynomials we will consider a rectangular grid $G\subseteq\mathbb{C}$ of the parameter space. Depending on the type of convergence of the orbit \eqref{orbitcirtic} we use different colors to give a graphic representation.
If the sequence $(R_\alpha^P)^k (q(y_\alpha))$ converges to an attracting $n$-cycle, then the sequence $L_{R^S}^{[r,s]}(\phi^S) (R_\alpha^S(y_\alpha))$ converges to a Lyapunov constant. We contemplate several cases:
\begin{enumerate}
\item  If $(R_\alpha^P)^k (q(y_\alpha))$ converges to one super-attracting fixed point
$[\frac{i}{2}:\frac{1}{2}],$ $[\frac{-i}{2}:\frac{1}{2}]$ or $[\frac{\alpha}{1 +\abs{\alpha}} : \frac{1}{1 +\abs{\alpha}}],$ then the point $\alpha $ of the parameter grid is colored  green, blue or denim, respectively. Note that the Lyapunov constant is equal to zero in this case.

\item If $(R_\alpha^P)^k (q(y_\alpha))$ converges to an $n$-cycle ($n>1$), then
$\alpha $ is colored red. In this case we may have attracting, indifferent or repelling $n$-cycles with Lyapunov constant being $<1$, $=1$ or $>1$, respectively.

\item If  $(R_\alpha^P)^k (q(y_\alpha))$ does not converge to any $n$-cycle, then the parameter value $\alpha $ is also colored red.
\end{enumerate}

\begin{figure}
\begin{center}
    \includegraphics[width=10cm]{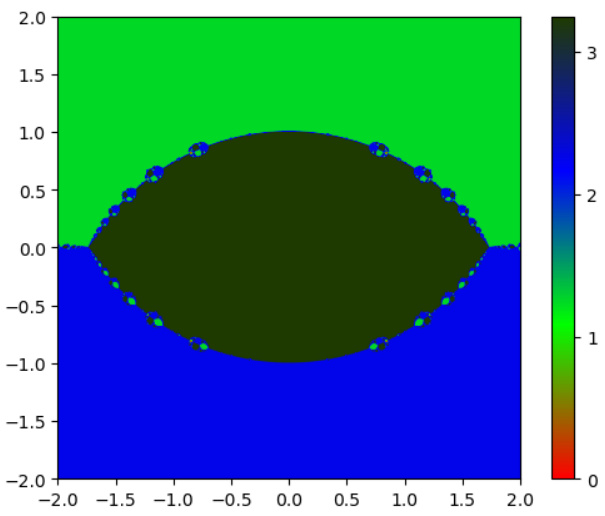}
    \caption{Parameter space of cubic polynomials associated to Newton's method}\label{cubicPolynomials1}
\end{center}
\end{figure}

\begin{figure}
\begin{center}
    \includegraphics[width=10cm]{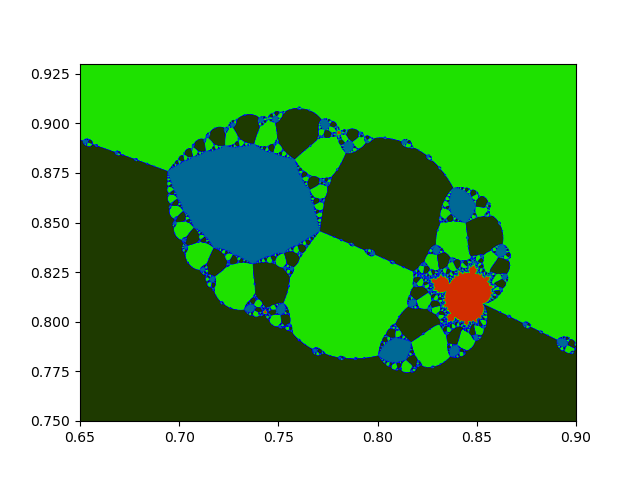}
    \caption{Fractal structure of the areas where a cyclic behaviour appears}\label{cubicPolynomials2}
\end{center}
\end{figure}

\begin{figure}
\begin{center}
    \includegraphics[width=10cm]{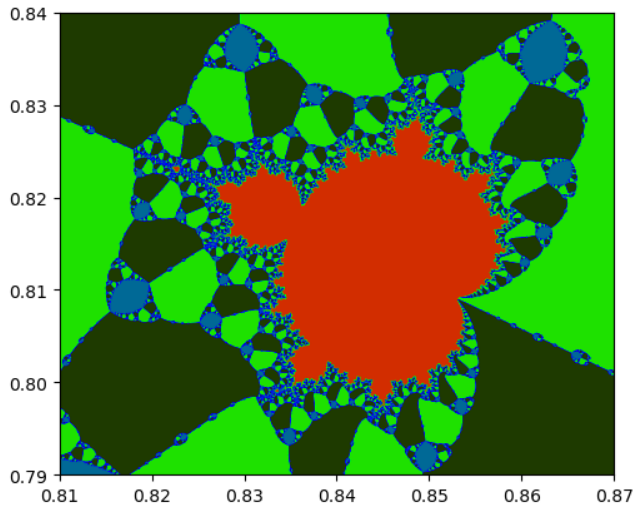}
    \caption{Mandelbrot-like set in which Newton's method fails}\label{cubicPolynomials3}
\end{center}
\end{figure}

Following the coloring strategy described above, one can compute graphical representations of any given region of $\mathbb{C}$ (we take rectangles in this case), as we show in Figures \ref{cubicPolynomials1}, \ref{cubicPolynomials2} and \ref{cubicPolynomials3}. Observe that the four colors previously detailed are also denoted by numbers in the graph color bar; red being $0$, green being $1$, blue being $2$ and blue denim being $3$.
In order to generate these graphics we have considered the maximum number of iterations of the Lyapunov-like technique (described in the previous sections) to be $50$, and a tolerance of $\frac{1}{10^8}$.

As we can see in Figure \ref{cubicPolynomials1}, at first sight, one cannot see any large red areas whose associated Newton's method behave ``badly''. It is also important to note that the region presented in
Figure \ref{cubicPolynomials1} is the rectangle $(-2,2)\times (-2,2)$.

Observe that the change of behaviour in Newton's method (i.e. change of color in the graphic) appears in an eye-shaped form in $(-2,2)\times (-2,2)$, and all along the real axis. Therefore, it is clear that, if any cyclic behaviour appears, it has to be where the behaviour of the polynomials changes.
Looking closer to one of the spots where the color changes, some visible red areas start to appear. In Figure \ref{cubicPolynomials2}, the rectangle $(0.65,0.9)\times (0.75,0.93)$ is shown. We have slightly reduced the blue color intensity in order to discern better the boundary around the red areas. In this figure, the fractal structure of the boundary areas is self-evident. Many small red areas can be also found.

In Figure \ref{cubicPolynomials3} (region $(0.81,0.87)\times (0.79,0.84)$) we can see the biggest red Mandelbrot-like area already observed in Figure \ref{cubicPolynomials2}. It should be noted that there are computational limitations, such as the maximum number of iterations as well as the given tolerance that we have had to choose. This is why the red area might not be a perfect-fit Mandelbrot-like set, when it conceptually should be. When zooming in other color-changing regions, one observes that every red area (that is, every area whose associated Newton's method behave ``badly'') has a similar fractal structure. We have also found some red areas where we know that the corresponding polynomial to each of its points presents an attracting $n$-cycle (or that the orbit of the arithmetic mean of its roots diverges) when Newton's method is applied.

One can still want to distinguish between each $n$-cycle by their period, in order to analyse the polynomials in which an attracting $n$-cyclic behaviour appears for a certain $n>1$. In this regard, we can modify the current graphics by adding some new colors. Again, computational limitations are necessary, so we have to choose the maximum period of the $n$-cycles that are going to be detected.
We may assign a different color for each $n$-cyclic behaviour. In this way, the area that corresponds to polynomials presenting a $2$-cyclic behaviour may be colored using one color, the area that corresponds to polynomials presenting a $3$-cyclic behaviour another different color, and so on for every $n$ between $2$ and the chosen maximum period.

\begin{figure}
\begin{center}
    \includegraphics[width=11cm]{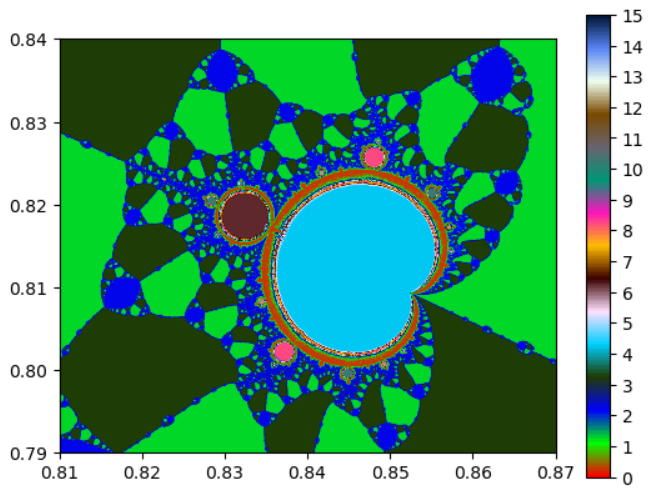}
    \caption{Distribution of  $n$-cyclic behaviour in a Mandelbrot-like set}\label{cubicPolynomials4}
\end{center}
\end{figure}

In Figure \ref{cubicPolynomials4}, colors $0$, $1$, $2$ and $3$ shown in the graph color bar are the same as in the previous figures. The cyan area (color $4$) visible in the main cardioid of the Mandelbrot-like set is composed by the points $\alpha $ for which the orbit of the critical point converges to an attracting $2$-cycle, the area in color $5$ by the points $\alpha $ for which the orbit of the critical point converges to an attracting $3$-cycle, and so on.
In  Figure \ref{cubicPolynomials4} we can see that each bulb of the Mandelbrot-like set has a clearly predominant area of an specific $n$-cyclic behaviour. For example, as we have mentioned, the main cardioid presents a $2$-cyclic behaviour predominance, whereas the bulb on its left has a predominant $4$-cyclic behaviour.
Also, the higher the maximum number of iterations allowed (that is, the higher the precision), the better these regions fit with $n$-cycles to each bulb. These graphic results suggest 
some theorical properties of these  bulbs.

\begin{remark}
It is not a coincidence that the areas whose associated polynomials present cyclic behaviour (i.e. where Newton's method fails to converge for some initial guess) are in fact Mandelbrot-like sets.
A more in-depth and rigorous study of this subject can be found in \cite{mcmullen_2000} and its references.
 Indeed,  in this work, C. T. McMullen proves  that every bifurcation set contains a copy of the boundary of a generalized Mandelbrot set of degree $d$, since generalized Mandelbrot sets are initial objects in the category of bifurcations. 
 
 In particular, the graphic results of our work confirm this fact, by considering the bifurcation locus of the holomorphic family of rational maps ($\alpha \in \C$) given by the application of Newton's method to cubic polynomials of the form $p(x)=(x-\alpha)(x^2+1)$, as previously mentioned.

\end{remark}

\subsection{Computing basins of n-cycles for the homogeneous Newton's method induced by the polynomial $p(x)=(x-0.845-0.81 i)(x^2+1)$}\label{example2basin}

We give a graphic representations of a rational map having a ``visible''  basin of one 2-cycle by using an implementation of the algorithms presented above. The rational map is associated to Newton's method and its implementation is given by normalized homogeneous pairs $R^S$ and $R^P$ given in equations \eqref{selfmaphopf} and \eqref{selfmaphopf1}. This means that, in order to obtain an implementation without any problems related to denominators close to zero,  we consider the homogeneous Newton's method $R={\rm HN}(p)$ as a self-map of the Hopf fibration.

In this particular case we have taken the cubic polynomial $p(x)=(x-0.845-0.81 i)(x^2+1)$. As a consequence of the study in the previous subsection, we have taken the value $\alpha=0.845+0.81 i$  by selecting a point of the main cardioid from the set ``similar'' to a Mandelbrot, see Figure \ref{cubicPolynomials3}. The values of $\alpha$ in the main cardioid correspond to rational maps having one attracting 2-cycle.
The implementation uses the Lyapunov functions $L_{\phi^S} ({\rm HN}(p))$, $L_{\phi^P} ({\rm HN}(p))$ where we take $\phi^S$, $\phi^P$  the spherical derivative of the rational function, which is computed by the method presented in subsection \ref{computingsphericalderivative} (see also  subsection  \ref{sphericalderivativesdiagram} and \eqref{cubical}).

In this example Newton's method associated to $\alpha=0.845+0.81\, i$ has as a free critical point $[\frac{\alpha}{3+\abs{\alpha}}: \frac{3}{3+\abs{\alpha}}] $ and its orbit converges to the attracting 2-cycle $$([ 0.289210 + 0.263312 \,i:1],[0.442483 + 0.986392\, i:1])=$$ 
\centerline{$(\small{[ \frac{0.289210 + 0.263312 \,i}{1.39112}:\frac{1}{1.39112}],[\frac{0.442483 + 0.986392\, i}{2.08109}:\frac{1}{2.08109}])}$}
\smallskip

\noindent which has Lyapunov constant equal to $ 0.501234$.
\begin{figure}
\begin{center}
    \includegraphics[width=10cm]{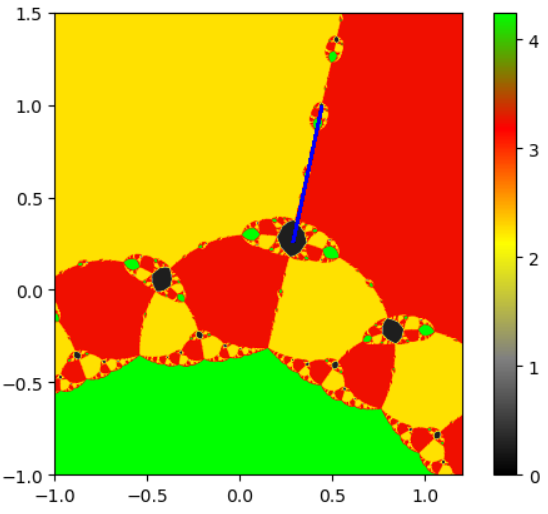}   \includegraphics[width=10cm]{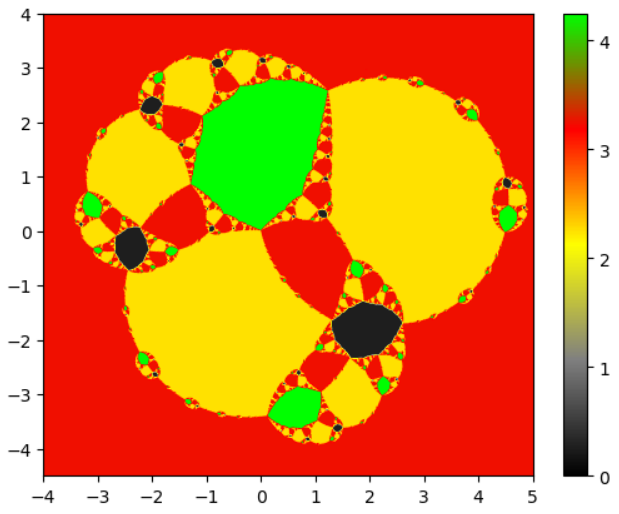}
    \caption{Basins of attraction of 3 fixed points and one 2-cycle}
    \label{basin2cycle}
\end{center}
\end{figure}

Figure \ref{basin2cycle} represents the basins of attraction  of three 1-cycles (fixed points) and one 2-cycle of $R^P$. At the top of the figure we have taken a rectangular neighbourhood at the origin $[0:1]$ and at the bottom of the figure we study the basins at a rectangular neighbourhood at the infinity point $[1:0]$. The legend contains the colors used to describe the basins:  Green color (4) stands for the basin of attraction of $[\frac{-i}{2}:\frac{1}{2}]$, yellow (2) for the basin of attraction of $[\frac{i}{2}:\frac{1}{2}]$, red (3) for the point $[0.845+0.81 \, i:1]=[\frac{0.845+0.81 \, i}{2.17052}:\frac{1}{2.17052}]$ and black (0) for the basin of an attraction  2-cycle, which is highlighted in the top image. The attracting 2-cycle stands out, represented as a cyan segment that joins the two poins of the 2-cycle.

We note that in Figure \ref{basin2cycle} the Julia set is not ``visible'' although it corresponds to the complement of the colored areas.

\section{Conclusions and future work}

In this work, we have presented a collection of algorithms and implementations based on the representation of a rational map as a rational self-map of the  Hopf fibration $S^{3} \to \p^{1}({\mathbb{C}})$ and on  the calculus of Lyapunov constants and exponents. Working with rational self-maps $(R^S,R^P)$ 
  allows us to compute the numerical  value of the function at any pole point and at the point at infinity. Moreover, one can avoid overflow, underflow and indetermination errors which usually appears computing  iterations of a rational map that is represented by an irreducible quotient of two univariate polynomials. 
  
The algorithms 
are  based on the calculus of Lyapunov constants and exponents associated to the spherical derivative of a rational map. This permits a nice connection with the ergodic theory of rational maps and the study of local stability of orbits of a rational map in the neighbourhood  at a point. We also give a canonical procedure to compute the spherical derivative which avoids the usual problems associated to the poles and the infinity point of the standard derivative.
 In order to avoid underflows which can appear computing limits of arithmetic means, we compute limits of the corresponding geometric means. We have used one implementation of these techniques in Julia language allows us a fast performance and permits numerical computations as well as  graphic constructions in a short period of time.
 We emphasize that the use of these techniques does not require  a previous calculus of fixed points of cycles of a given rational map.
 Finally, we give an example of the efficiency of these implemented algorithms  applying these methods to the study of the Newton's method for the case of cubic polynomials. 

Given a rational self-map $(R^S, R^P)$ of the Hopf fibration $S^3 \to \p(\C) $ one interesting question could be  the study of the set of  fixed points and cycles of the map $R^S \colon S^3 \to S^3$ and their  comparison with the  set of  fixed points and cycles of the map $R^P \colon \p(\C)  \to \p(\C) $. We also  think that these algorithms  can be modified and adapted  to give new  implementations  for computing  the Hausdorff dimension of Julia sets and other related numerical invariants. 
The techniques developed in this paper can be applied to other families of rational maps---for instance, the ones  obtained from 
Chebyshev's method and Halley's method (note that any iterative method represented by a rational map also admits a homogeneous version). The authors believe that more interesting future work can be developed by  applying these  implementations  in the study of bifurcation phenomena on parameter spaces.

\begin{remark}

We have prepared this work on the occasion of the 70th birthday of Professor Jos\'e Manuel Rodr\'{\i}guez Sanjurjo. Such work is framed in the area of discrete dynamics in which he and his collaborators have made important advances and contributions, see \cite{Sanjurjo20,Sanjurjo22} among many others. We want our humble contribution to serve as a tribute to Professor Sanjurjo and as an acknowledgment of a lifetime dedicated to teaching and research.

\end{remark}


\bibliography{Lyapunov_Hopffibration}
\bibliographystyle{elsarticle-harv}

\end{document}